*The restoration of Book X of the* Elements *to its original Theaetetean form*

by Stelios Negrepontis and Dimitrios Protopapas

**Abstract.** In the present work, we aim to restore Book X of the *Elements* to its original Theaetetean, pre-Eudoxean form in two separate ways.
First, we restore the considerable mathematical content of Book X, by correlating Book X with Plato's account of Theaetetus' mathematical discoveries and Plato's imitations of these discoveries for his philosophy. Thus, Theaetetus proved
(i) the eventual periodicity of the anthyphairesis of lines a to b, satisfying $Ma^2 = Nb^2$, for MN not square number, as deduced from Plato's *Theaetetus* and *Sophist*, and not simply their incommensurability with arithmetical means, as suggested by the mathematically flawed Proposition X.9;
(ii) the eventual periodic anthyphairesis of lines a to b, satisfying more general quadratic expressions, including the Application of Areas in defect, and employing this to show that 12 classes of alogoi lines, including the *minor*, despite being alogoi, are determined by an eventually periodic Application of Areas in defect, of relevance to the structure of the regular icosahedron in Book XIII of the *Elements*, and crucial for Plato's *Timaeus,* who has indicated his interest in the method in the *Meno* 86-87*;*
(iii) the anthyphairetic palindromic periodicity of the anthyphairesis of the surds √N for any non-square number N, as deduced from Plato's *Statesman*, not mentioned at all in Book X but containing all the essential mathematical tools for its proof, and of relevance to the general Pell's Diophantine problem, not mentioned in Book X but containing the essential mathematical tools for its proof.
Secondly, we restore the proofs of all propositions of Book X, in such way that these are proofs based on Theaetetus', and not on Eudoxus' theory of proportion of magnitudes, in particular not making any use of Eudoxus' condition (namely of definition 4 of Book V). The restoration is based on our reconstruction of Theaetetus' theory of proportion for magnitudes, for the limited class of ratios a/b such that either a, b are commensurable or the anthyphairesis of a to b is eventually periodic, without employing Eudoxus' condition, and its success provides a confirmation of our reconstruction.

**Section 1. Introduction**

In the present work, we aim to restore Book X of the *Elements* to its original Theaetetean, pre-Eudoxean form in two separate ways.

The first, outlined in Section 2, takes into account Plato's description of Theaetetus' epoch-making mathematical discoveries. Plato in several of his dialogues, including *Theaetetus, Sophist, Statesman, Meno*, and *Timaeus*, according to his own admission, imitated Theaetetus' mathematical discoveries on incommensurabilities in order to obtain the knowledge of the intelligible beings, the fundamental entities in Plato's philosophy. Plato's philosophical accounts of Theaetetus' contributions have been analyzed by Negrepontis, 2012, 2018, 2024, Negrepontis, Farmaki, Brokou, 2024, Negrepontis, Kalisperi, 2024, Negrepontis, Protopapas, 2025. It emerges from the analysis of Plato's account that Theaetetus was preoccupied with theorems about lines (namely straight-line segments) in periodic anthyphairesis, in particular it is suggested that he proved the following theorems:

> (1) If a, b are lines, M, N natural numbers, such that MN is not a square, and $Ma^2 = Nb^2$, then the anthyphairesis of a to b is eventually periodic, as suggested in the *Theaetetus* and the *Sophist* (Section 2.1).

There are no traces of theorem (1) in Book X of the *Elements*. Instead, Euclid presents

> Proposition X.9: If a, b are lines, M, N numbers, such that MN is not a square, and $Ma^2 = Nb^2$, then a, b are incommensurable,

a weak corollary of (1), with a faulty proof that, when corrected, appears to be related to Archytas' arithmetical methods in Book VIII, and not to Theaetetus' anthyphairetic methods.

> (2) If a, b are lines, N not a square number, and $a^2 = Nb^2$, then the anthyphairesis of a to b is eventually palindromically periodic,

as suggested in the *Statesman* (Section 2.2).

There are no traces of (2) in Book X, but it is remarkable that the tools for the proof of (2) are contained intact in Book X: the concepts of *apotome* lines (prime examples are all the remainders of the anthyphairesis of a to b in case $a^2 = Nb^2$ for non square N), of *binomial* lines, and the crucial *conjugacy* between them (Propositions X.112-114).

> (3) If ζ, η are lines commensurable in square only, ζ > η, and x is the solution of the Application of Areas in defect $x(ζ – x) = η^2/4$, then the anthyphairesis of ζ to x is eventually periodic (Section 2.3);

The 12 classes of alogoi lines, defined in Propositions X. 36-41 and X.73-78, are constructed much later, in Propositions X.54-59 and X.91-96, with the help of the solution x of the Application of Areas in defect (3); as a result, the 12 alogoi lines, despite their alogia, are nevertheless determined by an eventually periodic Application of Areas in defect, as hinted in the revealing *Scholia in Euclidem* X.135, 185.

The method of proof of (1) and (3) is basically Pythagorean, contained in the restored Book II of the *Elements (cf*, Negrepontis, Farmaki, 2024), plus a pigeon-hole argument

There are no traces of (3) in Book X, but (3) is employed to show that the alogos minor line, constructed in Proposition X.94, determines the structure of the canonical icosahedron (Proposition XIII.11), and is thus crucial for Plato's *Timaeus*, explaining Plato's interest in (3), as attested by the *Meno* 86e-87b passage.

Once (2) is established,

> (4) the solution of the general Pell problem (given a non square number N, to find natural numbers x and y such that $y^2 = Nx^2+1$)

is rather natural; in fact Book X contains evidence that Theaetetus was preoccupied with Pell's problem: one of the two criteria for the classification of apotome lines into six kinds is closely related to Pell's problem, the second is Proposition X.97, according to which if $\Omega$ is an apotome line, and $\Omega^2 = \gamma r$, then $\gamma$ is an apotome line, and the Pell number of $\gamma$ is the square of the Pell number of $\Omega$. Thus essentially Book X contains the tools for the proof of the Pell problem.

We regard Plato's account of Theaetetus' mathematical contributions as more authoritative than Euclid's, and thus we have no doubt that (1), (2), (3), (4) are theorems proved by Theaetetus. Thus, the first step towards the restoration of Book X to its original Theaetetean form is to add theorems (1), (2), (3), (4) on periodic anthyphairesis together with their proofs to Book X.

The reason why Euclid has eliminated all mention of periodic anthyphairesis is now not difficult to perceive: at the time Euclid wrote the *Elements*, Eudoxus' theory of proportion had prevailed over Theaetetus' theory; in fact periodic anthyphairesis, even though it continued to be the dominant concept in Plato's philosophy, had receded and was no longer an area of active study and research in Greek Mathematics. Euclid then, in composing Book X, followed the practice of Mathematics of his time, thus in approaching proofs of incommensurability he adopted the arithmetical methods initiated by Archytas, rather than the older method of periodic anthyphairesis that was dear to Theaetetus and Plato (Proposition X.9), and in questions of proportion of magnitudes, he adopted Eudoxus' arithmetized theory rather than Theaetetus' anthyphairetic theory of proportion for magnitudes, thus appearing not to side with those lamenting the prevalence of Eudoxus' theory as does the anonymous commentator in *Scholia in Euclidem* X.2. Euclid thus is seen, despite Proclus' claim of the contrary, not particularly of "Platonic persuasion".

In Section 3 we review the generally negative assessments of earlier scholars on Book X, assessments that are understandable, once it is seen that these assessments have not taken into essential consideration Plato's account of Theaetetus' contributions.

The second step of restoration rests on the realization that in Greek Mathematics there have been, not one but two theories of proportion for magnitudes, one based on equality of anthyphairesis, reported by Aristotle in *Topics* 158b, by all accounts by Theaetetus and another by Eudoxus. In Negrepontis, Protopapas, 2025, we have reconstructed Theaetetus' theory of ratios, only for the limited class of ratios a to b, such that the anthyphairesis of a to b is either finite or eventually periodic, and without any use of Eudoxus' condition. The basic properties of Theaetetus' theory of proportion are summarized in Section 4.

A necessary step toward restoring Book X to its original Theaetetean form, and at the same a confirmation of the correctness of our reconstruction of Theaetetus' theory of proportion is to show that all propositions of Book X of the *Elements* can be proved solely with the use of our reconstructed theory of proportion. This is achieved, in Section 5, by introducing modifications to the proofs in Euclid's *Elements*; in some cases these modifications are mild, in some other cases they are more essential (for instance the reader may compare our proofs of Propositions X.105 and X.113 with Euclid's), while for Propositions X.2, X.9 it has been necessary to devise altogether new proofs. The proofs of the propositions contained in what we regard as the final part of Book X (Propositions X.91-96 and X.54-59, in which we have a process of constructing an alogos line from an apotome/binomial line, and the converse process in Propositions X.97-102 and X.60-65 – related to the proof of Pell property – with the help of Propositions X.27-28, X.31-32 that lead to the application of areas in defect X.33-35) were the most demanding and delicate in their process and were thus the ones that required major alterations from their original, Euclidean form.

In Section 6, the restoration of Book X of the *Elements* to its original Theaetetean form according to these two criteria, of content and method of proof, and the resulting restoration of the correct historical order of the Books of the *Elements* are summarized in two tables.

## Section 2. Plato's account suggests that some fundamental propositions related to eventually periodic anthyphairesis and attributed to Theaetetus are absent from Book X, even though some of the tools needed for their proofs are there

Book X is the longest, and clearly the most complicated and most difficult to understand among all thirteen books of Euclid's *Elements*. Van der Waerden (1954, p. 172) correctly expresses the general sentiment:

> Book X does not make easy reading. The author succeeded admirably in hiding his line of thought.

Book X is generally attributed to Theaetetus by Plato in his trilogy *Theaetetus, Sophist, and Statesman,* by Pappus in his *Commentary on Book X of Euclid's Elements*, and by *Scholia in Euclidem* X.62, while *Scholia in Euclidem* XIII.1 implicitly correlates Theaetetus with the applications of Proposition X.94 to the structure the icosahedron in Proposition XIII.11. Thus, we would expect that Book X is a valid and authoritative source for Theaetetus' mathematical discoveries. However, we will argue that there are substantial differences between Plato's account of Theaetetus' achievements and the contents of Book X and that Plato's dialogues *Theaetetus, Sophist, Statesman*, where Plato, according to his admission, imitates philosophically Theaetetus' mathematical contributions, and also the *Meno* and the *Timaeus*, if properly understood, are a source far more authoritative and enlightening than Book X of the *Elements* and one able to reveal Theaetetus' line of thought and his contributions.

**2.1. Theaetetus' first theorem on the periodic anthyphairesis of lines a, b such that $Ma^2 = Nb^2$ with M, N natural numbers and MN a non-square number, as accounted philosophically by Plato in the *Theaetetus*, *Meno*, *Sophist***

**2.1.1.** *Some definitions and propositions*

To be able to describe Plato's account we must recall some basic definitions and propositions.

If a, b are two magnitudes, a > b, the *anthyphairesis* of a to b is the following process:

$a = k_0 b + c_1$, $b > c_1$

$b = k_1 c_1 + c_2$, $c_1 > c_2$

$c_1 = k_2 c_2 + c_3$, $c_2 > c_3$

…

$c_{n-1} = k_n c_n + c_{n+1}$, $c_n > c_{n+1}$

…

If some remainder $c_n$ divides the previous remainder $c_{n-1}$, then the anthyphairesis is *finite*, otherwise *infinite*; the anthyphairesis is (*eventually*) *periodic* if the sequence

$k_0, k_1, k_2, …, k_n, …$

of quotients is eventually periodic, namely if there are indices m < n, such that

$k_0, k_1, k_2, …, k_{m-1}$, period($k_m, …, k_n$);

and is (*eventually*) *palindromically periodic* if it is eventually periodic and the period is symmetric with respect to its middle.

Proposition X.2 of the *Elements* states that if the anthyphairesis of a to b is infinite, then a, b are *incommensurable*.

**2.1.2.** *Summary of Plato's philosophical account of Theaetetus' first theorem*
(a) Theaetetus' first theorem on quadratic incommensurabilities is stated in mathematically vague philosophical language (*Theaetetus* 147d-e);

(b) Socrates expressed the urge to imitate Theaetetus' mathematical discovery for the purpose of obtaining knowledge of the intelligible true Beings (*Theaetetus* 148d-5);
(c) Socrates hinted that the description of the intelligible Beings will have the form either "True Opinion plus Logos", or, equivalently, "Name and Logos" (*Theaetetus* 201d-202c);
(d) the full knowledge of the infinite anthyphairesis of the diameter to the side of a square is described as True Opinion/first two steps in the anthyphairesis, plus Logos/the Logos Criterion for anthyphairetic periodicity (*Meno* 80d-86a, 97a-98b*)*; and
(e) the full knowledge of the intelligible Beings "Angler" and "Sophist" are described as Name, a finite initial segment of philosophical anthyphairesis, plus Logos, the philosophical analogue of the Logos Criterion for anthyphairetic periodicity (*Sophist* 218b-221c, 264b-268d).

**2.1.3.** *The statement of Theaetetus' first theorem*

Line a is *commensurable* in square with line b if the squares $a^2$, $b^2$ with sides a, b respectively are commensurable (Definition X.I.2). In this case we may also say that line a is <u>rational</u> with respect to line b. Indeed, a line a is *rational* with respect to a given line r if $a^2$, $r^2$ are commensurable areas, otherwise the line a is *alogos/irrational* (Definition X.I.3).

Our analysis of Plato's account leads us to the conclusion that Theaetetus proved the following

> ***Proposition 2.1.*** *Theaetetus' first theorem.*
> If two lines a, b, with a > b, are such that $Ma^2 = Nb^2$, for natural numbers N, M, with NM non-square, then the anthyphairesis of a to b is eventually periodic, hence by Proposition X.2, a, b are incommensurable lines.

The method of proof that Theaetetus followed must have been a general form of the method followed by the Pythagoreans and by Theodorus, namely the propositions of Book II of the *Elements* restored to its initial Pythagorean form, heavy use of Pythagorean Application of Areas, including the preservation of the Gnomons, while the only really new idea of the proof must have been the pigeonhole principle (Negrepontis, Protopapas, 2025). The Pythagorean tools in a restored Book II of the *Elements* can be found in Negrepontis, Farmaki, 2024.

**2.1.4.** *Logos Criterion, the birth of Theaetetus' theory of proportion for magnitudes*

Even though Theaetetus' first theorem can be proved without any reference to a theory of ratios of magnitudes, we realize from Plato's account that Theaetetus introduced his theory of proportion of magnitudes, based, according to Aristotle's *Topics* 158b, on the *definition of proportion a/b = c/d* if the anthyphairesis of a to b is *equal* to the anthyphairesis of c to d, namely if the sequence Anth(a, b) of anthyphairetic quotients of a to b is equal to the sequence

Anth(c, d) of anthyphairetic quotients of c to d, mainly in order to express the awkward Pythagorean preservation of Gnomons (cf. 2.1.3, above), with the equivalent more flexible *Logos Criterion*, for the anthyphairetic periodicity. With the notation of Section 2.1.1, the *Logos Criterion* establishing periodic anthyphairesis can be expressed as follows:

> there are two indices m < n such that $c_m/c_{m+1} = c_n/c_{n+1}$.

**2.1.5.** *Proposition X.9 of the* Elements *juxtaposed to Theaetetus' first theorem*

By contrast, Euclid in the *Elements* presents Proposition X.9, stating that if two lines a, b, with a > b, are such that $Ma^2 = Nb^2$, for natural numbers N, M, with NM non-square, then a and b are incommensurable.
We observe that although the assumptions of Proposition X.9 are similar to those of the first theorem of Theaetetus, the conclusion is much weaker.
The proof of Proposition X.9, as Mazur, 2007 has pointed out, is mathematically unsound; when corrected, it is seen to employ Archytas' weaker arithmetical methods, presented in Book VIII of the *Elements.* (Negrepontis, 2018; Negrepontis, Farmaki, Brokou, 2024; Negrepontis Farmaki, Kalisperi, 2022; Negrepontis, Protopapas, 2025, Sections 4, 5, 7, 9).
It is remarkable that Proposition X.2, although stated and proved in the *Elements* (with a proof however that makes unnecessary use of Eudoxus' condition) is not used for the proof of Proposition X.9, or in fact of any other proof, and practically has no role in the *Elements*! We must assume that it was left over from a previous version of Book X, in which anthyphairesis had a role.
Even though most historians of Greek Mathematics have adopted variants of the proof of Propostion X.9 (cf. Negrepontis, Protopapas, 2025, Section 9), there is no question that Plato's account of Theaetetus' first mathematical discovery is more authoritative than the one by Euclid in the *Elements*.

**2.2. Theaetetus' second theorem on palindromically periodic anthyphairesis of surds, described as a philosophical imitation in Plato's *Statesman*, is fully absent from Book X**

**2.2.1**. *Plato's imitation of anthyphairetic palindromic periodicity in the* Statesman

The whole of Plato's *Statesman,* the sequel of the *Sophist* in the trilogy *Theaetetus, Sophist, Statesman,* is devoted to the knowledge of the intelligible Being "Statesman", by a method that is a direct continuation of the method Name and Logos, followed in the *Sophist*, and thus clearly an imitation, a philosophical analogue of periodic anthyphairesis. The *palindromic periodicity theorem* is described philosophically in the *Statesman* by a division that, being the direct continuation of the divisions in the *Sophist*, is clearly anthyphairetic. It is clear that Plato's philosophical imitation can only be the imitation of the following

> ***Proposition 2.2.*** *Theaetetus' second theorem,*
> If the lines a, b, are such that $a^2 = Nb^2$, with N non-square number
> (namely they form a surd),
> then the anthyphairesis of line a to line b is palindromically periodic.

This is a fundamental theorem proved in modern times by Lagrange (1769, 1771) and Euler (1765). But how could Theaetetus have proved this theorem, a fact strongly suggested by Plato's *Statesman*?

**2.2.2.** *We realize that the tools for the proof of the palindromic periodicity are in Book X! Apotome lines, binomial lines, and their conjugation*

A line γ is an *apotome* with respect to a given line r if there are two lines ζ, η, with ζ > η and such that ζ, η are incommensurable but $ζ^2$, $η^2$, $r^2$ are commensurable and γ = ζ – η (definition in Proposition X.73).
In the *Elements* it is nowhere explained why these lines are defined, but we find that

> all the *remainders* $c_n$ of the anthyphairesis of a to b, in case $a^2 = Nb^2$ for some non-square number N, are apotome lines with respect to b.

The crucial step in Book X is the introduction of the binomial line, the corresponding additive line to the subtractive apotome: A line δ is a *binomial* with respect to a given line r if there are two lines ζ, η, with ζ > η and such that ζ, η are incommensurable but $ζ^2$, $η^2$, $r^2$ are commensurable and δ = ζ + η (definition in Proposition X.36).
The significance of the binomial lines is that they serve as the *conjugate* lines to the apotome lines (Propositions X.112-114) in the following exact sense:
Now it is shown that

> both an apotome and a binomial line are *alogoi* lines with respect to r
> (Propositions X.73, 36, respectively),
> BUT the product of an apotome line γ = ζ – η with the corresponding conjugate
> binomial line δ = ζ + η produces a rational area $ζ^2 – η^2$ (Propositions X.112-114).

This powerful quadratic conjugation served as the inspiration to Bombelli (1572) for conceiving the complex numbers, and as noted by Weil (1984), it provides the principal tool for the proof of the *palindromic periodicity theorem*. Thus, its role in Book X is to provide the tools for the proof of the palindromic periodicity of the anthyphairesis/continued fraction of the surd √N for every non-square number N (Negrepontis, 2018; Negrepontis, Farmaki, Brokou, 2024).

**2.2.3.** It is remarkable that Euclid, although providing the mathematical tools for its proof, is hiding completely Theaetetus' original intention, an intention that can be revealed and restored only by taking into account our understanding of Plato's philosophical imitation in the *Statesman*.

**2.3. The 12 alogoi lines in Book X of the *Elements* are, despite their alogia, held together, controlled by periodic anthyphairesis**

**2.3.1.** *From the apotome lines to the subtractive and from the binomial lines to the additive alogoi lines*
The 12 *alogoi* lines consist of

6 subtractive, depending on 6 kinds of apotome lines, defined in Propositions X.73-78, and 6 additive, depending on 6 corresponding kinds of binomial lines defined in Propositions X.36-41.

For each kind of apotome or binomial line γ with respect to the given line r, an *alogos line* Ω of that kind is determined as the mean proportional of γ and r: $Ω^2 = γr$ (X.91-96 for apotomes producing subtractive alogoi lines, X.54-59 for binomials producing additive alogoi lines).

**2.3.2.** *The role of the Application of Areas in defect in establishing the relation between the 6 kinds of apotome and the 6 subtractive alogoi lines, and of the 6 kinds of binomial lines and the 6 additive alogoi lines.*

The study of the alogoi lines in Book X depends heavily on the Pythagorean Application of Areas in defect, essentially contained in Book II of the *Elements*, but Euclid has removed from Book II the Pythagorean Proposition II.5/6 because there is the more general Proposition VI.28, involving Eudoxus' theory of ratios of magnitudes.

Given an apotome line γ = ζ – η, we construct the Application of Areas in defect $x(ζ – x) = η^2/4$, and we solve this problem geometrically by constructing in the circle of diameter ζ = AB a line of length η/2 = DE at right angles to the diameter and ending in the circumference of the circle at the point D. We form the right triangle ADB with perpendicular sides Φ = AD, Ψ = BD and hypotenuse ζ = AB. Then $Φ^2 = AB·AE = ζ(ζ – x)$, $Ψ^2 = AB·BE = ζx$, and Ω = Φ – Ψ.

**2.3.3.** *Theaetetus' first theorem generalized: If γ = ζ – η is an apotome, then the anthyphairesis of ζ, x, satisfying the Application of Areas in defect equation $x(ζ – x) = η^2/4$, is eventually periodic*

We observe crucially that by Theaetetus' first theorem,

> the anthyphairesis of ζ, x is eventually periodic,

Indeed, with the same method of proof as for Theaetetus' first theorem (Section 2.1.3), the following theorem about the Pythagorean Application of Areas in Defect can be proved.

> If two lines ζ, η, ζ > η, are commensurable in square only and x is the solution of the Application of Areas in defect equation $x(ζ – x) = η^2/4$, then the anthyphairesis of ζ to x is eventually periodic.

The general form of Theaetetus' first theorem is examined in Negrepontis, Protopapas, 2025, Section 6.

**2.3.4**. *Scholion in Euclidem X.135 emphasizes that the 12 alogoi lines, despite their alogia, are controlled and dominated by eventually periodic anthyphairesis*

Since $Φ^2 = ζ(ζ – x)$, $Ψ^2 = ζx$, and Ω = Φ – Ψ, it is clear that Ω is determined by x, and x is constructed by an eventually periodic anthyphairesis. This then makes understandable the *Scholion in Euclidem X.135, 1-9*

> It is worthy of admiration (Θαυμάζειν ἄξιον), that

    the *staying/holding together power(κρατητικὴ δύναμις)*
    of the triad (τῆς τριάδος)
    *finitizes (ἀφορίζει) the alogos power* (ἄλογον δύναμιν)
    and reaches (διήκει) till the extremes (μέχρι τῶν ἐσχάτων),
    because each kind of *alogia* (ἀλογίας) is
    *finitized* (ἀφορίζεται) by some mean (μεσότητος),
    one by the geometric (γεωμετρικῆς) mean,
    the other by the arithmetical (ἀριθμητικῆς) mean,
    the other by the harmonic/musical (μουσικῆς) mean;
    and it seems that the Being of the soul (ἡ τῆς ψυχῆς οὐσία)
    sittting upon (ἐπιβατεύουσα) (the substance of) the magnitudes
    according to the *"logoi"* of the soul (κατὰ τοὺς ἐν αὐτῇ λόγους) *finitizes* (ὁρίζειν)
    everything indefinite (ἀόριστον) and the infinite of the alogia (τὴν τῆς ἀλογίας ἀπειρίαν)
    by means of these three *bonds* (δεσμοῖς) [translation by the authors]

Indeed the "staying power", "sitting upon the magnitudes", "finitizing alogia" in *Scholion* X.135 is precisely this eventually periodic anthyphairesis (cf. Negrepontis, Farmaki, Kalisperi, 2022).

Thus, although an alogos line Ω is a "bad" line, being alogos, nevertheless the twelve alogoi lines, generally denoted by Ω, are held together, controlled, dominated by an eventually periodic process.

**2.3.5.** The most interesting case is the construction of the *minor* line Ω from the *fourth apotome* line γ (Proposition X.94). This is essentially the only Proposition of Book X that has an application outside Book X, namely in the construction of the icosahedron in Proposition XIII.11. Indeed, the side of the icosahedron is a "minor line" with respect to the diameter of the sphere in which the icosahedron is comprehended, defined in X.76 and constructed in Proposition X.94.

**2.3.6.** *Plato's interest for the construction of the alogoi lines by means of an Application of Areas in defect is revealed by* Meno *86e-87b*

Negrepontis and Farmaki (2019, pp. 288-296) initially correlated passage 86e-87b of Plato's *Meno* to Proposition X.33 and eventually to Proposition X.94, while Negrepontis, Farmaki and Kalisperi (2022, Section 8.3) argued, on linguistic grounds, that the triangle described in the passage must be a *right angled triangle* and not an isosceles one, as usually depicted. Thus, *Meno*'s passage indicates Plato's interest for the part of Book X that employs Application of Areas in defect and leads to Proposition XIII.11 on the structure of the icosahedron.

Plato's interest in Book XIII is understandable, since the regular polyhedra play a fundamental role in the *Timaeus* (cf. Lamb, 1925; Negrepontis, Kalisperi, 2024). The binding element here is the fact that some key three-dimensional space constructions such as the construction of the regular icosahedron in Book XIII have some connections with periodical anthyphairesis and are in this sense finitized, cf. *Scholia in Euclidem* X.135, 185.

**2.4. The solution of the general Pell Diophantine equation $y^2 = Nx^2 + 1$ for every non-square number N is suggested by the palindromically periodic anthyphairesis of every**

**surd (1.2), the classification of the apotome and the binomial lines, and Proposition X.97 in Book X**

**2.4.1.** The *classification of an apotome (or a binomial) line* $\gamma = \zeta - \eta$ ($\delta = \zeta + \eta$, correspondingly), with respect to a given line r, in six kinds, given in X.84/85 (X.47/48 for the binomials), depends on two criteria:
*criterion (i)* on whether r is commensurable, either to $\zeta$, or to $\eta$, or to neither of the two, and
*criterion (ii)* on whether the line $\theta$, defined by the rule $\theta^2 = \zeta^2 - \eta^2$, is or is not commensurable to the line $\zeta$.

**2.4.2.** Euclid does not explain the role of these criteria, and it does not seem that modern scholars have provided any explanation, but our analysis has shown that both of these criteria are relevant to the apotome lines generated as the anthyphairetic remainders of the anthyphairesis of a to b for $a^2 = Nb^2$, with N a non square number, with respect to the given line b; indeed, then
*criterion (i)* depends on the anthyphairetic interchange between odd and even stages of the anthyphairesis, but more importantly
*criterion (ii)* is closely related to *the Pell number* of the apotome line.

**2.4.3.** Proposition X.97 implies that if $\Omega$ is the alogos line apotome, and $\Omega^2 = \gamma r$, then the Pell number of $\gamma$ is the square number of the Pell number of $\Omega$. With exactly the same method, it can be proved that if $\Omega_1, \Omega_2$ are apotome lines and $\Omega_1\Omega_2 = \gamma r$, then the Pell number of $\gamma$ is the product of the Pell numbers of $\Omega_1, \Omega_2$.

**2.4.4.** A solution of the general Pell problem can be obtained employing the palindromic periodicity theorem (2.2, above), criterion (ii) for the classification of apotome (X.84/85), and binomial lines (2.4.1, above), and Proposition X.97 and Remarks (2.4.3, above) are relevant to the solution of the Pell problem for every non-square N (Negrepontis, Farmaki, Brokou, 2024).

**2.5.** We now have a satisfactory picture of the original mathematical contributions by Theaetetus, based (a) on our interpretation of Plato's dialogues, and (b) on the considerable mathematical tools that we find in Book X of the *Elements,* notably the conjugation of apotome and binomial lines, and the Pell related classification of apotome lines.

Theaetetus proved essentially three theorems on periodic anthyphairesis (a first periodicity theorem 2.1, the palindromic periodicity for the special case of surds 2.2 and a generalized periodicity theorem 2.3.3), and used the third one (2.3.3) in one direction for his approach to the general Pell problem (2.4); while again in the converse direction for proving a fundamental property of the icosahedron in Proposition XIII.11 (2.3).

It must be noted that all these mathematical discoveries of Theaetetus are propositions about lines in periodic anthyphairesis, and their proofs employ Theaetetus' theory of proportion defined only for ratios, the terms of which are in finite or eventually periodic anthyphairesis.

This picture differs radically from Book X itself, where every trace of periodic anthyphairesis has been obliterated. Propositions 2.1, 2.2, and the implications of Section 2.4.3 do not appear in Book X, only Proposition X.9, an incommensurability result receiving a (faulty but) strictly arithmetical Archytan-type proof, while all proofs of the propositions in Book X are now based solely on Eudoxus' theory of proportion of magnitudes. Nevertheless the essential mathematical tools described in (b) with which Theaetetus proved these theorems, notably the conjugation of apotome and binomial lines, and the Pell related classification of apotome lines, are kept in Book X.

Thus Book X has powerful mathematical tools, which can be used for the proof of fundamental propositions closely related to periodic anthyphairesis, but the propositions themselves are missing and nowhere stated; in addition the proofs of these tools employ Eudoxus' theory of ratios and not a theory of ratios based on periodic anthyphairesis, as was the case with Theaetetus' theory, further distancing these mathematical tools from the propositions they were meant to prove. The negative comments on Book X expressed by earlier scholars (presented in section 3) are then quite understandable, granted that these scholars never related Book X of the *Elements* with Plato's philosophy.

Euclid in the *Elements* has obliterated any mention of Theaetetus' theory of proportion of magnitudes, treating all proportion of magnitudes in terms of Eudoxus' theory of ratios (Books V and VI of the *Elements*); in this he was surely following the common mathematical practice of his time: Eudoxus' theory had prevailed over Theaetetus' theory. However, Book X, authored by Theaetetus, was certainly based on his own theory of proportion; thus at some point Euclid himself, or someone before Euclid, possibly Hermotimus of Colophon as Knorr conjectures (1983, p. 59), modified all proofs of the propositions of Book X, so as to conform to Eudoxus' theory of proportion.

But we must have in mind that Eudoxus' theory of ratios was a component of the program, conceived by Archytas and Eudoxus, to take Mathematics away from periodic anthyphairesis, right after it reached its apex in the Mathematics of Theaetetus and in the philosophy of Plato. It became, if not rigorously proved, intuitively clear that periodicity in anthyphairesis was tied to plane geometry and quadratic equations, and not to space geometry and cubic equations, and that periodic anthyphairesis was not suitable to deal with cubic incommensurabilities, or with a problem like the duplication of the cube. Thus, the mathematician that modified the original work by Theaetetus in Book X of the *Elements*, replaced proofs based on propositions of Theatetetus' theory of ratios of magnitudes with proofs based on propositions of Eudoxus' theory of proportion, a rather easy task since Eudoxus' theory is far stronger than Theaetetus' theory; but, going even further, that mathematician evidently considered propositions about periodic anthyphairesis and the Logos Criterion itself as irrelevant (even though one could prove the Logos Criterion for anthyphairetic periodicity within Eudoxus' theory, e.g. using Proposition VI.16 (stating that for any lines a, b, c, d, we have a/b = c/d if and only if ad = bc), and thus replaced the two theorems of Theaetetus on periodic anthyphairesis 2.1 and 2.2, with the far weaker Proposition X.9 on incommensurability, by means of an arithmetical proof.

As a result of these modifications, Book X was transformed from the original work of Theaetetus, dealing with theorems on periodic and on palindromically periodic anthyphairesis, employing a theory of ratios of magnitudes based on periodic anthyphairesis and the Logos

Criterion, and the conjugacy of apotome and binomial lines, and approaching Pell's problem, with the theorem of palindromic periodicity and the Pell related classification into six apotome lines, with statements of propositions and their proofs that had no trace of periodic anthyphairesis, but keeping intact the mathematical tools that were leading to the proofs of these propositions.

Thus, our goal to restore Book X to the original Theaetetean work from which it evolved consists of two different tasks. First, to include an account of all Theaetetus' discoveries and not only those that have an application to the structure of the icosahedron (Proposition X.94). Secondly, to revert all the proofs of the propositions in Book X to proofs employing solely Theaetetus' theory of proportion for magnitudes, with no use whatsoever of Eudoxus' principle. This will be accomplished in Section 5, below.

At this point it is necessary to comment on Proclus, *Scholia in Euclidem* 68,20-23 where he writes:

> καὶ τῇ προαιρέσει δὲ Πλατωνικός ἐστι καὶ τῇ φιλοσοφίᾳ ταύτῃ οἰκεῖος, ὅθεν δὴ καὶ τῆς συμπάσης στοιχειώσεως τέλος προεστήσατο τὴν τῶν καλουμένων Πλατωνικῶν σχημάτων σύστασιν.

> Euclid belonged to the persuasion of Plato and was at home in this philosophy; and this is why he thought the goal of the *Elements* as a whole to be the construction of the so-called Platonic figures. [translation by Morrow, 1970]

Euclid's presentation of Book X of the *Elements* shows that Euclid had no real philosophical affinity to Plato's philosophy.

## Section 3. The negative comments on Book X of the *Elements* by earlier scholars

The differences between Plato's account of Theaetetus' contributions and Euclid's Book X are thus striking. We realize that if someone studies only Book X, without gaining a true understanding of Plato's account, then it is nearly impossible to arrive at a true assessment of Theaetetus' epoch making mathematical discoveries. These differences explain the negative comments expressed by earlier scholars on Book X, since in fact they had no satisfactory interpretation of Plato's account.

Indeed, Book X has been studied by earlier scholars and historians of Greek Mathematics, starting with Bombelli and Stevin in the 16th century, and their assessments are generally negative. Despite the early genial association by Bombelli (1572) of the conjugation of complex numbers with the conjugation of the apotome and binomial lines in Propositions X.112-114 of the *Elements*, an association that was crucial for grasping the concept of complex numbers, there is a long list of scholars, including Stevin (1585), Heiberg (1883-1885), Heath (1926), van der Waerden (1954), Mueller (1981), Taisbak (1982), Knorr (1983), and even Weil (1984), that have assessed Book X as a work of great difficulty ("the cross of mathematicians"), but without goal, point, utility, or purpose, unless it is "hidden admirably", and even as a work having no mathematical importance, concealing no mathematical treasures under its obscurity, and in fact "a mathematical blind alley", a veritable "cul de sac".

**3.1. The great idea by Bombelli on the essential connection between the conjugacy of the apotome and binomial lines in Book X of the *Elements* and the conjugacy of complex numbers**

Weil (1984, p. 202) describes beautifully the connection between the basic concepts of Book X, the apotome (introduced in X.73) and the binomial line (introduced in X.36), and their conjugacy in Propositions X.112–114, a conjugacy very similar to the relation between a complex number and its conjugate, that impressed Bombelli in 1572.

> The first concept, motivated by the work of the Italian algebraists of the sixteenth century on equations of degree 3 and 4, had been introduced by Bombelli in Book I of his Algebra of 1572, in close imitation of Euclid's theory of irrationals of the form $a + \sqrt{b}$ in his Book X ('binomial' in the terminology of Campanus' Latin translation of Euclid; accordingly, binomio is Bombelli's word for a complex number $a + \sqrt{-b}$). Further developments belong to the history of algebra and analysis rather than to number theory; suffice it to say that Euler played a decisive role in extending to complex numbers the main operations of analysis.

**3.2. The negative opinions on Book X of the *Elements* by Stevin, 1585, Heiberg, 1883-1885, Heath, 1926, van der Waerden, 1954, Mueller, 1981, Taisbak, 1982, Knorr, 1983, and even Weil, 1984**

Stevin, 1585, likens Book X with the "cross of the mathematician" because of its difficulty, although it has no utility. Heath, 1926, seems to believe that the main use of Book X is in the study of the five regular polyhedra in Book XIII. Van der Waerden, 1954, finds Book X not easy to read, and declares that its author "succeeded admirably in hiding his thought". Mueller, 1981 finds that Book X has "no clear mathematical goal" and in fact constitutes "a mathematical blind alley". Taisbak, 1982, believes that Book X "has no other point than to entertain us with good logic", "has no mathematical importance". Knorr, 1983, believes that Book X contains no hidden "mathematical treasures" and writes: "The student who approaches Euclid's Book X in the hope that its length and obscurity conceal mathematical treasures is likely to be disappointed", agreeing with the characterization assigned to Book X by Mueller and Taisbak as "a blind alley".

Especially for Propositions X.112-114, Heiberg, 1883-1885, regards them as later interpolations (Heath and Knorr disagree), while Weil and Knorr express their disappointment over Euclid's utter failure to perceive the potential of these propositions.

Stevin, 1585

> The difficulty of the Tenth Book of Euclid has become for many in horror, even to calling it the cross of the mathematicians, a subject matter too hard to digest and in which they perceive no utility. [translation from French by Knorr, 1983, p.61]

Heiberg, 1883-1885, in his classical edition of Euclid's *Elements,* considers Propositions X.112-114, as later interpolations. Heath,1926, p. 246, writes:

> Heiberg considers that this proposition and the succeeding ones are interpolated, though the interpolation must have taken place before Theon's time. His argument is that X.112-115 are nowhere used, but that X.111 rounds off the complete discussion of the 13 irrationals (as indicated in the recapitulation), thereby giving what was necessary for use in connexion with the investigation of the five regular solids.

Heath, 1926, p.9

> It will naturally be asked, what use did the Greek geometers actually make of the theory of irrationals developed at such length in Book X.? The answer is that Euclid himself, in Book XIII, makes considerable use of the second portion of Book X. dealing with the irrationals affected with a negative sign, the apotomes etc.
> …
> Of course the investigation in Book X. would not have been complete if it had dealt only with the irrationals affected with a negative sign. Those affected with the positive sign, the binomials etc., had also to be discussed, and we find both portions of Book X., with its nomenclature, made use of by Pappus in two propositions, of which it may be of interest to give the enunciations here.

Van der Waerden, 1954, p.172

> Book X does not make easy reading ... The author succeeded admirably in hiding his line of thought.

Mueller, 1981, p.270-271

> One would, of course, prefer an explanation that invoked a clear mathematical goal intelligible to us in terms of our own notions of mathematics and which, under analysis, would lead univocally to the reasoning in book X. Unfortunately, book X has never been explicated successfully in this way nor does it appear amenable to explication of this sort. Rather, book X appears to be an expedient for dealing with a particular problem and at the same time a mathematical blind alley.

Taisbak, 1982

> The output of the lesson will be a logic game with coloured quadrangles and line segments, leading to an unimportant classification of some line segments. Apart from the fun one can have (and some Greek surely had) from manipulating such objects, we shall not hesitate to maintain that the game, bar very few details, has no mathematical importance; the X'th book of

> the *Elements* may well be called a *cul-de-sac* in mathematics, even though it did inspire Kepler to his *Harmoniae Mundi*.
> All the same, it is a fascinating, nay haunting, piece of literature, and its composition reminds one of epic or dramatic poetry, with its long sequences of uniform statements apt to be learnt by heart and transmitted orally. p. 27.

> We are prepared to face the possibility that there was no other point than to entertain us with good logic p. 58

> We cannot but applaud van der Waerden's judgement, that "the author succeeded admirably in hiding his line of thought". It may have been part
> of the teaching: if one knows (as we do) what it is all about, the order
> of the *Elements* is as easy to follow as any didactic and pedagogic presentation,
> perhaps even easier to learn by heart; and whoever is not initiated may be more easily kept
> outside by the opacity of the theme. A peculiar instance of Μηδείς Εισίτω. p.66

Knorr, 1983

> It is clear that Euclid's closing theorem [[X.115]] suddenly throws open the field of inquiry, after the main body of the theory had been neatly tied and sealed in the scholium to X, 111 only four theorems earlier. This may be good cause for viewing the result as a post-Euclidean addition, as the prominent editors Heiberg and Heath do. [Heath [1926, III, 255]]. Heath resists Heiberg's similar suspicions about the product theorems in 112-114, and I believe rightly so [Ibid., p. 246]. p.56

Knorr then argues that the completed theory of Book X is neither by Theaetetus nor by Eudoxus, but by some minor mathematician in the period between Eudoxus and Euclid, in fact probably by Hermotimus of Colophon (1983, p.58-59).

> one modern judgment of the work [[Mueller [1981, 271]] as "a mathematical blind alley" is quite apt.
> Certain potentials of the subject matter are entirely missed, for instance, insight into how products of the form $(a \pm b\sqrt{N})(c \pm d\sqrt{N})$ bear on the finding of solutions to integral relations of the form $x^2 - Ny^2 = \pm m$.
> [[At the time of his course of lectures on the history of the Pell equation
> (Institute for Advanced Study, 1978-79), Professor A. Weil expressed to me
> his disappointment over Euclid's utter failure to perceive this potentially
> fruitful development of the theory of Book X. The meager hints we might glean
> from Archimedes and Diophantus are little encouragement for supposing that
> the Greeks recognized these possibilities or pursued them to any length.]] p.60

> The student who approaches Euclid's Book X in the hope that its length and obscurity conceal mathematical treasures is likely to be disappointed. As we have seen, the mathematical ideas are few and capable of far more perspicuous exposition than is given them here. The true merit of Book X, and I believe it is no small one, lies in its being a unique specimen of a fully elaborated deductive system of the sort that the ancient philosophies of mathematics consistently prized. p.60

Weil, 1984, p. 15, realizes the potential importance of the presence of conjugacy in Book X, but at the end dismisses the possibility that these conjugacies were actually used in the manner that they inspired Bombelli:

> Not only is Euclid himself well aware of the relation $(\sqrt{r} + \sqrt{s})(\sqrt{r} - \sqrt{s}) = r - s$ but even the identity $1/(\sqrt{s} + \sqrt{r}) = \sqrt{r}/(r - s) - \sqrt{s}/(r - s)$ may be regarded as the essential content of prop. 112 of that book. Unfortunately, Euclid's motivation in Book X seems to have been the wish to construct a general framework for the theory of regular polygons and polyhedra, and not, as modern mathematicians would have it, an algebraic theory of quadratic fields. So we are left to speculate idly whether, in antiquity or later, identities involving square roots may not have been used, at least heuristically, in arithmetical work.

## Section 4. Theaetetus' theory of proportion for magnitudes is for a limited enough class of ratios, and has no need of Eudoxus' condition

We now turn to the second step of the restoration of Book X and we posit that *the first part of Book X in its original form was Theaetetus' theory of ratios of magnitudes*. This theory has been reconstructed in Negrepontis and Protopapas, 2025, and it must be kept in mind that our reconstruction of the theory does not employ in any of the proofs Eudoxus' principle and for the sake of completeness, we include here the basic propositions of the reconstruction. It is clear that Euclid has suppressed all reference to Theaetetus' theory in favor of Eudoxus' theory (correspondingly Books V and VI of the *Elements*). Again, this step is towards obtaining the original form of Book X, in place of the one modified by Euclid.

We note that there is no trace of Theaetetus' theory of ratios of magnitudes in Book X, in fact in all of the *Elements*, and no trace of periodic anthyphairesis. It is clear that Euclid has already introduced deep changes in the content of Book X, following the Archytas-Eudoxus motion away from Theaetetean periodic anthyphairesis, and hence away from Plato's philosophy.

The detailed reconstruction of the proof of Theaetetus' theorem reveals that Theaetetus' theory of magnitudes is *essentially a theory of ratios of lines*, designed principally for the needs of periodic anthyphairesis. Critical for the development of Theaetetus' theory of ratios was the introduction by Theaetetus of the "Logos Criterion" (Negrepontis, 2012; Negrepontis, Protopapas, 2025, Section 2).

**Proposition 4.1.1** (*Logos Criterion for the periodicity of an anthyphairesis*). Let a, b be two line segments, with $a > b$, and anthyphairesis given as follows:
$a = k_0 b + e_1$, with $b > e_1$,
$b = k_1 e_1 + e_2$, with $e_1 > e_2$,
…

$e_{m-1} = k_m e_m + e_{m+1}$, with $e_m > e_{m+1}$,
$e_m = k_{m+1} e_{m+1} + e_{m+2}$, with $e_{m+1} > e_{m+2}$,
…
$e_{n-1} = k_n e_n + e_{n+1}$, with $e_n > e_{n+1}$,
$e_n = k_{n+1} e_{n+1} + e_{n+2}$, with $e_{n+1} > e_{n+2}$,
…

and assume that there are indices m < n such that $e_n/e_{n+1} = e_m/e_{m+1}$ ("Logos Criterion"). Then the anthyphairesis of a to b is eventually periodic, in fact Anth(a, b) = [$k_0$, $k_1$, …, $k_m$, period ($k_{m+1}$, $k_{m+2}$, …, $k_n$)].

The Logos Criterion for anthyphairetic periodicity is essentially just a restatement in more suggestive and convenient language of the more heavy and burdensome Pythagorean preservation of Application of Areas/Gnomons. As shown in Negrepontis, Protopapas, 2025, Section 8, Theaetetus' theory holds for the limited class of ratios of lines with finite or eventually periodic anthyphairesis.

Exactly because of the limited nature of the theory, it is possible to prove, with the help of the crucial Proposition 4.3.4, all the basic properties of this theory with no recourse to Eudoxus' condition (as was the case with the previous reconstructions), and this leads to the establishment of the Theaetetean analogue 4.3.5 of Proposition V.9 of the *Elements*, without any use of Eudoxus' condition.

The extension of the theory from lines to areas and surfaces is achieved with the help of the Theaetetean analogue 4.3.11 of the key Proposition VI.1 of the *Elements*. The fact that this theory is really a theory not for general magnitudes, but principally for lines, makes necessary the many special cases in the proof of its Propositions. This feature is pointed out by Aristotle in *Analytics Posterior* 74a, and the vehicle of passing from lines to areas and surfaces is exactly the proposition (analogue to Proposition VI.1 of the *Elements*) singled out by Aristotle in the *Topics* 158b to inform us about this theory.

## 4.1 The transition from the Pythagorean preservation of Gnomons to the equality of the anthyphaireses

The propositions of this section depend on

**Proposition** (*Elements* I.44)**.** Given lines a, b, c, to construct a line d, such that ab = cd.

**Proposition 4.1.2.** If a, b, c, d are lines and ad = bc, then Anth(a, b) = Anth(c, d).

*Proof.* Let Anth(a, b) = [$k_0$, $k_1$, …]. We proceed by induction.
Thus $a = k_0 b + e_1$, with $e_1 < b$. Then $cb = ad = k_0 bd + e_1 d$.
By Proposition I.44 of the *Elements*, there is a line segment $f_1$, such that $e_1 d = b f_1$.

Since $e_1 < b$, it is clear that $f_1 < d$.

Thus $cb = k_0 bd + e_1 d = k_0 bd + bf_1$; then $c = k_0 d + f_1$ with $f_1 < d$.

Thus $\text{Anth}(a, b) = [k_0, k_1, \ldots] = [k_0, \text{Anth}(b, e_1)]$, while $\text{Anth}(c, d) = [k_0, \text{Anth}(d, f_1)]$.

Thus, we have $e_1 d = bf_1$ and we must prove that $\text{Anth}(b, e_1) = \text{Anth}(d, f_1)$.

We continue as in the first step, and we finish the proof by induction.

**Lemma 4.1.3.** If $Aa^2 = Bab + Cb^2$ and $Aa^2 = Bad + Cd^2$, then $b = d$.

*Proof.* By hypothesis $Bab + Cb^2 = Bad + Cd^2$. Suppose b is different from d, say $b > d$.
Then $Bab > Bad$ and $Cb^2 > Cd^2$, hence $Bab + Cb^2 > Bad + Cd^2$, a contradiction.

**Proposition 4.1.4.** If $Aa^2 = Bab + Cb^2$ and $Ac^2 = Bcd + Cd^2$, then $ad = bc$.

The condition of Proposition 4.1.4 is referred to as *the preservation of Application of Areas/Gnomons*.

*Proof.* By Proposition I.44, there is a line e such that $bc = ae$. From the assumption we have
$Aa^2 c = Babc + Cb^2 c$, hence
$Aa^2 c = Baae + Cbae$, hence
$Aac = Bae + Cbe$, hence
$Aacc = Baec + Cbec$, hence
$Ac^2 a = Bcea + Ce^2 a$, hence
$Ac^2 = Bce + Ce^2$.
By Lemma 4.1.3, $d = e$. Hence $ad = bc$.

**Proposition 4.1.5.** If $Aa^2 = Bab + Cb^2$ and $Ac^2 = Bcd + Cd^2$, then $\text{Anth}(a, b) = \text{Anth}(c, d)$.

*Proof.* If $Aa^2 = Bab + Cb^2$, and $Ac^2 = Bcd + Cd^2$, then $ad = bc$ (Proposition 4.1.4); and if $ad = bc$, then $\text{Anth}(a, b) = \text{Anth}(c, d)$ (Proposition 4.1.2).

### 4.2. The introduction of the generalized side and diameter numbers

#### 4.2.1. The generalized side-diameter numbers $p_n$, $q_n$

**Definition 4.2.1.** Let there be an infinite anthyphairesis of a to b:
$a = k_0 b + e_1$, $e_1 < b$
…
$e_{m-2} = k_{m-1} e_{m-1} + e_m$, $e_m < e_{m-1}$
$e_{m-1} = k_m e_m + e_{m+1}$, $e_m < e_{m+1}$
$e_m = k_{m+1} e_{m+1} + e_{m+2}$, $e_{m+1} < e_{m+2}$
…
$e_{n-2} = k_{n-1} e_{n-1} + e_n$, $e_n < e_{n-1}$
$e_{n-1} = k_n e_n + e_{n+1}$, $e_{n+1} < e_n$

… .

We set $p_1 = 1$, $q_1 = k_0$, and $p_n = k_{n-1}p_{n-1} + p_{n-2}$, $q_n = k_{n-1}q_{n-1} + q_{n-2}$ for every n.

*Note.* The (classical) side and diameter numbers are identified as the generalized side and diameter numbers of the anthyphairesis of the diagonal (diameter) to side of a square.
[Indeed the anthyphairesis of diameter to side of a square is $k_0 = 1$, $k_n = 2$ for all $n = 1, 2, …$ .
Then the generalized side-diameter numbers are recursively defined by $p_1 = 1$, $q_1 = k_0 = 1$, and,
$p_{n+2} = 2p_{n+1} + p_n$, $q_{n+2} = 2q_{n+1} + q_n$.
Then $p_{n+2} = 2p_{n+1} + p_n = p_{n+1} + p_{n+1} + p_n = p_{n+1} + (p_n + q_n) + p_n = p_{n+1} + q_{n+1}$, and
$q_{n+2} = 2q_{n+1} + q_n = q_{n+1} + q_{n+1} + q_n = q_{n+1} + 2p_n + q_n + q_n = q_{n+1} + 2(p_n + q_n) = 2p_{n+1} + q_{n+1}$, and
the recursive definition of the side and diameter numbers is derived
$p_{n+2} = p_{n+1} + q_{n+1}$, $q_{n+2} = 2p_{n+1} + q_{n+1}$].

### 4.2.2. The generalized side and diameter numbers are useful in expressing anthyphairetic remainders in terms of the initial magnitudes a and b

**Proposition 4.2.2** (the anthyphairetic remainders in terms of the generalized side and diameter numbers). With the notation of Definition 4.2.1, we have $e_n = (-1)^n(q_n b - p_n a)$ for all $n = 1, 2, …$

*Proof.* $e_1 = a - k_0 b = (-1)^1(q_1 b - p_1 a)$.
Assume inductively that $e_{n-2} = (-1)^{n-2}(q_{n-2}b - p_{n-2}a)$, and $e_{n-1} = (-1)^{n-1}(q_{n-1}b - p_{n-1}a)$.
Then $e_n = e_{n-2} - k_{n-1}e_{n-1} =$
$(-1)^{n-2}(q_{n-2}b - p_{n-2}a) - (-1)^{n-1}k_{n-1}(q_{n-1}b - p_{n-1}a) =$
$(-1)^n[q_{n-2}b - p_{n-2}a + k_{n-1}(q_{n-1}b - p_{n-1}a)] =$
$(-1)^n[(q_{n-2} + k_{n-1}q_{n-1})b - (p_{n-2} + k_{n-1}p_{n-1})a] =$
$(-1)^n(q_n b - p_n a)$.

**Proposition 4.2.3.** If $Aa^2 = Bab + Cb^2$, and $ad = bc$, then $Ac^2 = Bcd + Cd^2$.

*Proof.* From the assumption we have $Aa^2 = Bab + Cb^2$,
hence $Aa^2c = Babc + Cb^2c$,
hence $Aa^2c = Babc + Cbad$,
hence $Aac = Bbc + Cbd$,
hence $Aac^2 = Bbc^2 + Cbdc$,
hence $Aac^2 = Bcad + Cad^2$,
hence $Ac^2 = Bcd + Cd^2$.

**Proposition 4.2.4.** (1) For every finite ordered sequence of natural numbers $k_0, k_1, ..., k_n$, there are lines a and b, $a > b$, such that $Anth(a, b) = [period(k_0, k_1, …, k_n)]$, and $p_{n+1}a^2 = (q_{n+1} - p_n)ab + q_n b^2$.
(2) If lines c and d satisfy $Anth(c, d) = [period(k_0, k_1, …, k_n)]$, then $p_{n+1}c^2 = (q_{n+1} - p_n)cd + q_n d^2$.

*Proof.* (1) The first $n + 1$ steps of the anthyphairesis of a to b are:

$a = k_0 b + e_1, e_1 < b$

$b = k_1 e_1 + e_2, e_2 < e_1$

$e_1 = k_2 e_2 + e_3, e_3 < e_2$

…

$e_{n-2} = k_{n-1} e_{n-1} + e_n, e_n < e_{n-1}$

$e_{n-1} = k_n e_n + e_{n+1}, e_{n+1} < e_n$.

We also impose the cross product condition $ae_{n+1} = be_n$. The cross product condition implies the periodicity of the anthyphairesis; indeed, by Proposition 4.1.2, Anth(a, b) = Anth($e_n$, $e_{n+1}$), and thus by induction: Anth(a, b) = [$k_0, k_1, \ldots, k_n$, Anth($e_n, e_{n+1}$)] =

[$k_0, k_1, \ldots, k_n$, Anth(a, b)] =

[$k_0, k_1, \ldots, k_n, k_0, k_1, \ldots, k_n$, Anth($e_n, e_{n+1}$)] =

… =

[period($k_0, k_1, \ldots, k_n$)].

The cross product relation further insures that a, b satisfy an Application of Areas in excess; indeed, by Proposition 4.2.2, $e_n = (-1)^n(q_n b - p_n a)$, and $e_{n+1} = (-1)^{n+1}(q_{n+1} b - p_{n+1} a)$, hence, $ae_{n+1} = (-1)^{n+1} a(q_{n+1} b - p_{n+1} a)$, and $be_n = (-1)^n b(q_n b - p_n a) = (-1)^n b(p_n a - q_n b)$.

Thus, the cross product condition takes the form: $a(q_{n+1} b - p_{n+1} a) = b(p_n a - q_n b)$, hence $p_{n+1} a^2 = (q_{n+1} - p_n) ab + q_n b^2$.

(2) Since Anth(a, b) = Anth(c, d), the anthyphairesis of c to d is purely periodic. Since the generalized side and diameter numbers depend solely on the sequence of quotients $k_0, k_1, \ldots, k_n$, it follows that the generalized side $r_i$ and diameter $s_i$ numbers of c to d coincide with the generalized side $p_i$ and diameter $q_i$ numbers of a to b, respectively, namely $r_i = p_i$, $s_i = q_i$ for i = 1, 2, …, n, n + 1; hence, denoting by $f_i$ the ith anthyphairetic remainder of c to d, it follows that $cf_{n+1} = df_n$, hence, by Proposition 4.2.3, c, d satisfy the same Application of Areas in excess, therefore $p_{n+1} c^2 = (q_{n+1} - p_n) cd + q_n d^2$. From this result, by completing the anthyphairesis of both a to b and c to d upwards, we obtain, without difficulty, the general case. The details are left to the reader.

### 4.3. Theaetetus' theory of proportion for magnitudes

Proposition 4.2.4 makes possible the proof of the fundamental Proposition 4.3.4: in Theaetetus' theory, for lines a, b, c, d, the condition of proportion a/b = c/d is equivalent to the equality of the cross products ad = bc, with no recourse to Eudoxus' condition. In turn, following essentially an old idea of Becker (1933), we can prove not only the *Alternando* property for lines (if a/b = c/d, then a/c = b/d), but also the crucial analogue of Proposition V.8 of the *Elements*.

**Definition 4.3.1.** A line a possesses Theaetetean ratio with respect to line b if the anthyphairesis of a to b is finite or eventually periodic.

**Definition 4.3.2** (of proportion)**.** Two Theaetetean ratios of lines a/b and c/d are equal in the Theatetean sense, a/b = c/d, if Anth(a, b) = Anth(c, d).

**Proposition 4.3.3** (analogue of the Transitive property, proposition V.11 for lines; the analogous proposition for numbers is not explicitly stated in Book VII, but is implicitly used)**.** Let a, b and c, d and e, f be pairs of lines possessing Theaetetean ratios, namely with finite or eventually periodic anthyphairesis. Then if a/b = c/d, and c/d = e/f, then a/b = e/f.

*Proof.* By Definition 4.3.1, either Anth(a, b) = [$k_0, k_1, ..., k_n$, period($k_{n+1}, k_{n+2}, ..., k_m$)] if it is eventually periodic, or Anth(a, b) = [$k_0, k_1, ..., k_n$] if it is finite. Also, by Definition 4.3.2, a/b = c/d means that Anth(a, b) = Anth(c, d). By Definition 4.3.2 again, c/d = e/f means that Anth(c, d) = Anth(e, f). Thus Anth(a, b) = Anth(e, f) and this, by the same definition, means that a/b = e/f.

**Proposition 4.3.4** (a fundamental proposition of our reconstruction and the analogue of Proposition VI.16 for lines, Proposition VII.19 for numbers)**.** Let a, b and c, d be pairs of lines possessing Theaetetean ratios, namely with either finite or eventually periodic anthyphairesis. Then a/b = c/d (namely Anth(a, b) = Anth(c, d)) if and only if ad = bc.

*Proof.* Let a, b and c, d be two pairs of lines with either finite or eventually periodic anthyphairesis. If ad = bc, then by Proposition 4.1.2, Anth(a, b) = Anth(c, d), namely a/b = c/d (Definition 4.3.2). If a/b = c/d, then by Proposition 4.2.4, a, b and c, d satisfy the same Application of Areas in Excess, say $Aa^2 = Bab + Cb^2$ and $Ac^2 = Bcd + Cd^2$. By Proposition 4.1.4, ad = bc.

**Proposition 4.3.5** (analogue of Proposition V.9 for lines; the analogous proposition for numbers is not explicitly stated in Book VII of the *Elements*, but is implicitly used)**.** If lines a, b and a, c have finite or eventually periodic anthyphairesis, and a/b = a/c, then b = c.

*Proof.* By Definition 4.3.2, a/b = a/c means that Anth(a, b) = Anth(a, c). By Proposition 4.3.4, ac = ab, hence b = c.

*Note.* For the restrictive class of lines a, b and c, d with finite or eventually periodic anthyphairesis, the proofs of the Fundamental Proposition 4.3.4 and 4.3.5 are provided without any use of Eudoxus' condition, Definition V.4 of the *Elements*, thus avoiding the anachronism that plagues the earlier reconstructions. This by itself is a strong argument in favor of our reconstruction.

Regarding the next proposition we will prove, the *Alternando*, van der Waerden writes:

> But one property, viz. the interchange of the means, causes difficulty. And now it is curious that, according to Aristotle, it was indeed exactly the proof of this proposition, which at first led to difficulties. "Formerly", says Aristotle in Anal. Post. 15, "this proposition was proved separately for numbers, for line segments, for solids and for periods of time. But after the introduction of the general concept which includes numbers as well as lines, solids and periods of time" (viz. the concept of magnitude), "the proposition could be proved in general".

The proof for numbers can be found in Book VII (VII 13). What was the old proof for line segments?

O. Becker has advanced an ingenious hypothesis for this. From the proportionality

(1) a : b = c : d,

one deduces first the equality of the areas

(2) ad = bc,

then interchanges b and c, and finally returns to the proportionality

(3) a : c = b : d. [1954, p. 177].

Van der Waerden refers to the following remarkable Aristotle's *Analytics Posterior* 74a17-25 passage:

> καὶ τὸ ἀνάλογον ὅτι καὶ ἐναλλάξ,
> ᾗ ἀριθμοὶ καὶ ᾗ γραμμαὶ καὶ ᾗ στερεὰ καὶ ᾗ χρόνοι,
> ὥσπερ ἐδείκνυτό ποτε χωρίς,
> ἐνδεχόμενόν γε κατὰ πάντων μιᾷ ἀποδείξει δειχθῆναι·
> ἀλλὰ διὰ τὸ μὴ εἶναι ὠνομασμένον τι ταῦτα πάντα ἕν,
> ἀριθμοί μήκη χρόνοι στερεά,
> καὶ εἴδει διαφέρειν ἀλλήλων,
> χωρὶς ἐλαμβάνετο. νῦν δὲ καθόλου δείκνυται·
> οὐ γὰρ ᾗ γραμμαὶ ἢ ᾗ ἀριθμοὶ ὑπῆρχεν,
> ἀλλ' ᾗ τοδί, ὃ καθόλου ὑποτίθενται ὑπάρχειν.

> Again, the law that proportionals alternate might be supposed to apply
> to numbers qua numbers, and similarly to lines, solids and periods of time;
> as indeed it used to be demonstrated of these subjects separately.
> It could, of course, have been proved of them all by a single demonstration,
> but since there was no single term to denote the common quality of
> numbers, lengths, time and solids,
> and they differ in species from one another, they were treated separately;
> but now the law is proved universally;
> for the property did not belong to them
> qua lines or qua numbers,
> but qua possessing this special quality which they are assumed to possess universally.
> [translation by Tredennick (1960)]

We will follow Becker's idea for the proof of

**Proposition 4.3.6** (*Alternando*, analogue of Proposition V.16 for lines, Proposition VII.13 for numbers)**.** If a/b = c/d and Anth(a, c), Anth(b, d) are finite or eventually periodic, then a/c = b/d.

*Proof.* By Proposition 4.3.4, ad = bc. Since Anth(a, c) is finite or eventually periodic, then by Proposition 4.3.4, a/c = b/d.

**Proposition 4.3.7** (*Ex Equali,* analogue of Proposition V.22 for lines, Proposition VII.14 for numbers). If a, b, c, d, e, f are lines and a/b = d/e and b/c = e/f and Anth(a, c), Anth(d, f) are finite or eventually periodic, then a/c = d/f.

*Proof.* By Proposition 4.3.4, ae = bd and bf = ce. Then abf = ace = bcd, hence af = cd. By Proposition 4.3.4, a/c = d/f.

**Proposition 4.3.8** (*perturbed proportion*, analogue of Proposition V.23 for lines). If a, b, c, d, e, f are lines and a/b = e/f and b/c = d/e, and Anth(a, c), Anth(d, f) are finite or eventually periodic, then a/c = d/f.

*Proof.* By Proposition 4.3.4, af = be and be = cd. By transitivity (Proposition 4.3.3), af = cd. By Proposition 4.3.4, we have a/c = d/f.

**Proposition 4.3.9** (part (a) is the analogue of Proposition V.12 for lines and Proposition VII.12 for numbers while part (b) is the analogue of Proposition V.19 for lines and Proposition VII.11 for numbers). If a, b, c, d are lines and Anth(a, b) = Anth(c, d) is finite or eventually periodic anthyphairesis, then:

(a) Anth(a + c, b + d) is finite or eventually periodic, and (a + c)/(b + d) = a/b, and

(b) if a > c and b > d, then Anth(a − c, b − d) is finite or eventually periodic, and (a − c)/(b − d) = a/b.

*Proof.* (a) By Proposition 4.3.4, and since Anth(a, b) = Anth(c, d), we have: ad = bc.
Then: ad + ab = bc + ab, a(b + d) = b(a + c). By Proposition 4.1.2, Anth(a, b) = Anth(a + c, b + d) and then, by Definition 4.3.2, (a + c)/(b + d) = a/b.
(b) By Proposition 4.3.4, and since Anth(a, b) = Anth(c, d), we have ad = bc.
Then ab + ad = ab + bc, ab − bc = ab − ad, b(a − c) = a(b − d). By Proposition 4.1.2, Anth(a, b) = Anth(a − c, b − d) and then, by Definition 4.3.2 (a − c)/(b − d) = a/b.

**Proposition 4.3.10** (part (a) is the analogue of Proposition V.18 for lines and part (b) is the analogue of Proposition V.17 for lines, but there is no analogue for numbers). If a, b, c, d are lines and Anth(a, b) = Anth(c, d) is finite or eventually periodic, then:

(a) Anth(a + b, b) = Anth(c + d, d) is finite or eventually periodic, and (a + b)/b = (c + d)/d, and

(b) if a − b > b and c − d > d, then Anth(a − b, b) = Anth(c − d, d) is finite or eventually periodic, and (a − b)/b = (c − d)/d.

*Proof.* (a) Consider the case where Anth(a, b) = Anth(c, d) and both are finite. Let their anthyphairesis be denoted by [$k_0, k_1, \ldots, k_n$]. Then, it is easy to demonstrate that the anthyphairesis of a + b to b and of c + d to d differs from the anthyphairesis of a to b and c to d, respectively, only in the first quotient. Specifically, instead of $k_0$, the first quotient now becomes $k_0 + 1$.
Then Anth(a + b, b) = Anth(c + d, d) = [$k_0 + 1, k_1, \ldots, k_n$] and, by Definition 4.3.2,
(a + b)/b = (c + d)/d. Similarly, we can demonstrate that the proposition is true when the anthyphairesis of a to b is eventually periodic.

(b) Consider the case where Anth(a, b) = Anth(c, d) and both are finite. Let their anthyphairesis be denoted by $[k_0, k_1, \ldots, k_n]$. Then, it is easy to demonstrate that the anthyphairesis of a − b to b and of c − d to d differs from the anthyphairesis of a to b and c to d, respectively, only in the first quotient. Specifically, instead of $k_0$, the first quotient now becomes $k_0 − 1$.

Then Anth(a − b, b) = Anth(c − d, d) = $[k_0 − 1, k_1, \ldots, k_n]$ and, by Definition 4.3.2, (a − b)/b = (c − d)/d. Similarly, we can demonstrate that the proposition is true when the anthyphairesis of a to b is eventually periodic.

**Proposition 4.3.11** (*Topics Proposition*, analogue of VI.1 for lines, Propositions VII.17, 18 for numbers). If a, b, c are lines, with Anth(a, b) eventually periodic, then (ac)/(bc) = a/b.

*Proof.* It is not difficult to prove that Anth(a, b) = Anth(ac, bc).
Indeed, let the anthyphairesis of a to b be the following:
$a = k_0 b + c_1, c_1 < b$,
$b = k_1 c_1 + c_2, c_2 < c_1$,
$c_1 = k_2 c_2 + c_3, c_3 < c_2$,
…
$c_{n-1} = k_n c_n + c_{n+1}, c_{n+1} < c_n$,
… .
Then
$ac = k_0 bc + c_1 c, c_1 < b$,
$bc = k_1 c_1 c + c_2 c, c_2 < c_1$,
$c_1 c = k_2 c_2 c + c_3 c, c_3 < c_2$,
…
$c_{n-1} c = k_n c_n c + c_{n+1} c, c_{n+1} < c_n$,
… .
It is clear that Anth(a, b) = Anth(ac, bc), hence by Definition 4.3.1, (ac)/(bc) = a/b.

**Proposition 4.3.12** (analogue of Proposition V.9 for rectilinear areas). If A, B, C are rectilinear areas, with A/B = A/C, and Anth(A, C) is finite or eventually periodic, then B = C.

*Proof.* By the definition of A/B = A/C, we have Anth(A, B) = Anth(A, C). Let r be a given line. By Propositions I.44-I.45 of the *Elements*, there are lines a, b, c, such that A = ar, B = br, C = cr. By Proposition 4.3.11, A/B = (ar)/(br) = a/b, A/C = (ar)/(cr) = a/c.
Hence, by transitivity (Proposition 4.3.3), a/b = a/c and by Proposition 4.3.4, ac = ab, hence b = c.
Hence B = br = cr = C.

**Proposition 4.3.13** (*Alternando,* analogue of Proposition V.16 for rectilinear areas). If A, B, C, D are rectilinear areas, with A/B = C/D, and Anth(A, C), Anth(B, D) are finite or eventually periodic, then A/C = B/D.

*Proof.* By the definition of A/B = C/D, we have Anth(A, B) = Anth(C, D). Let r be a given line. By Propositions I.44-I.45 of the *Elements*, there are lines a, b, c, d, such that A = ar, B = br, C = cr, D = dr. By Proposition 4.3.11, A/B = (ar)/(br) = a/b, C/D = (cr)/(dr) = c/d.

Hence, by transitivity (Proposition 4.3.3), a/b = c/d and by Proposition 4.3.4, ad = bc. By Proposition 4.3.6 (*Alternando* for lines) a/c = b/d and by Proposition 4.3.11 and transitivity (Proposition 4.3.3), A/C = B/D.

**Proposition 4.3.14** (Analogue of Proposition V.22, part (a) is the *Ex Equali* for rectilinear areas and part (b) is the *Ex Equali* for a mixed proportion). (a) If A, B, C, D, E, F are rectilinear areas, with A/B = D/E, and B/C = E/F, and if Anth(A, C), Anth(D, F) are finite or eventually periodic, then A/C = D/F.
(b) If A, B, C are rectilinear areas, d, e, f are lines, with A/B = d/e, B/C = e/f, and if Anth(A, C), Anth(d, f) are finite or eventually periodic, then A/C = d/f.

**Proposition 4.3.15** (Analogue of Proposition V.23, part (a) is the *perturbed proportion* for rectilinear areas and part (b) is the *perturbed proportion* for a mixed proportion). (a) If A, B, C, D, E, F are rectilinear areas, with A/B = E/F, and B/C = D/E, and Anth(A, C), Anth(D, F) are finite or eventually periodic, then A/C = D/F.
(b) If A, B, C are rectilinear areas, d, e, f are lines, with A/B = e/f, B/C = d/e, and Anth(A, C), Anth(d, f) are finite or eventually periodic, then A/C = d/f.

The proofs of Propositions 4.3.14, 4.3.15 are similar to the proof of Proposition 4.3.12 and are left to the reader.

*Note*. The fragmentation of the proof of a property in separate cases for lines, areas, volumes, as described by Aristotle in *the Analytics Posterior* passage quoted above, is confirmed, not only for the *Alternando* property, but also for the other properties of Theaetetus' theory of proportion.

## Section 5. Theaetetus' theory of proportion for magnitudes is sufficient for Book X of the *Elements*

We recall our comments in Negrepontis, Protopapas (2025, Section 3) on the reconstructions of Theaetetus' theory of proportion for magnitudes by Becker (1933), van der Waerden (1954) and Knorr (1975):

> Becker, 1933 assumes that Theaetetus' theory of ratios of magnitudes, based on the definition of proportion in terms of equal anthyphairesis, is a theory for the same class of magnitudes as is Eudoxus' theory. Therefore, his reconstruction of Theaetetus' proof of the analogue of the crucial Proposition V.9 of the Elements, is inevitably in need of Eudoxus' condition (Definition V.4); van der Waerden, 1954, p.178-179, essentially agrees with this reconstruction.
> …

Knorr, 1975, p.340, although not fully accepting Becker's reconstruction, nevertheless does not question Becker's assumption that Theaetetus' theory is about the same class of magnitudes as is the theory of Eudoxus, and that therefore a rigorous proof of the analogue of Proposition V.9 would need Eudoxus' condition.

In particular, according to these scholars, Book X of the *Elements* depends on Eudoxus' condition (Definition 4 in Book V of the *Elements*) since the proofs in Book X make use of Eudoxus' theory of proportion for magnitudes (Books V and VI of the *Elements*).

Mueller (1981, p.292) explicitly states that the foundations of Book X of the *Elements* include Eudoxus' Books V and VI of the *Elements*:

> The important foundations of book X, other than definitions, may be characterized as follows: (i) Standard proportion theory (the laws of book V) and its geometric applications in VI, 1, 8, 16, 17 (14 in X,22), 19-20, 22.

On the contrary, Taisbak (1982) has the "haunting impression" that Book X does not need Eudoxus' theory of proportion, but states that "such an impression cannot be proved, perhaps not even made plausible":

> But it is a haunting impression (which haunted Oskar Becker already in the 1930's) that the Story of the Greater and the Lesser Line, viz. the X'th book, did not need Eudoxus' theory of proportion in order to be coherently told. Such an impression cannot be proved, perhaps not even made plausible; but we entertain the thesis that it is possible to tell that story within a theory of ratio that is confined to commensurable magnitudes, supplied by very few (less than 5) "intuitively true" statements about the ratio of line segments and of quadrangles, whether they be commensurable or no. [p. 17-18]
>
> The theorem [[X.11]] is problematic in a pre-Eudoxean context (cf. § 1.2), since we have not defined proportionality of magnitudes that are not commensurable; but we shall use it of such magnitudes only in connection with
> T VI 1 Rectangles under the same height are to each other as their bases.
> …
> Proposition VI 1 is crucial to the theory of proportion in the set of plane polygons.
> We demand that it be granted us, without proof, in the rectangular form (which may well be the one that Aristotle alludes to in Topica 158 b 29 ff.)
> Also without proof we accept some "intuitive" propositions, which were undoubtedly believed to hold for incommensurable as well as for commensurable magnitudes, even before Eudoxus;
> T V 7 Equal quadrangles have the same ratio to one and the same quadrangle. And one quadrangle has the same ratio to equal quadrangles.
> T V 9 Quadrangles which have the same ratio to one and the same quadrangle, are equal.
> V 11 Ratios which are equal to one and the same ratio, are equal. [p. 22-23]

The litmus test for the validity of our reconstruction in Negrepontis, Protopapas, 2025, of Theaetetus' theory of ratios of magnitudes is whether or not our reconstructed theory is sufficient for the proof of all the propositions of Book X of the *Elements*. This is so because Book X is Theaetetean, and because Euclid, as van der Waerden points out, has replaced the Theaetetean theory of ratios that Theaetetus must certainly have used, with the Eudoxean theory of Books V and VI of the *Elements*. If our reconstruction is valid, we must be able to go through all the propositions of Book X using only the limited, compared to the Eudoxean, Theaetetean theory of ratios, namely the theory applying
*only to the limited class of magnitudes a, b, such that Anth(a, b) is finite or eventually periodic.*

## 5.1. Theaetetus' theory of ratios for magnitudes

We believe that all the propositions from Theaetetus' theory of ratios of magnitudes must have originally been a part of the content of Book X in its Theaetetean form. Specifically, we believe that we must begin by placing to this position, the very first part of a restored Book X, our reconstruction of the theory of ratios developed by Theaetetus, exactly as we have documented it in Section 4.

## 5.2. Propositions X.1-9, related to Theaetetus' first incommensurability theorem

### 5.2.1. The statement of Proposition X.2

If the anthyphairesis of magnitude a to magnitude b is infinite, then a, b are incommensurable.

### 5.2.2. The Euclidean proof of Proposition X.2

The proof employs Proposition X.1, whose proof is based on Eudoxus' condition. Thus, it would appear that there would be a problem in employing this proposition in connection with Theaetetus' incommensurabilities.

### 5.2.3. A proof of Proposition X.2 with no use of Eudoxus' condition

However, there is a quite natural and elementary proof of Proposition X.2 using the fundamental Pythagorean Propositions VII.1 & 2 of the *Elements* to find the greatest common divisor of two numbers by the method of anthyphairesis.
*Propositions VII.1 & 2.*

If a > b are two unequal natural numbers, then (i) the anthyphairesis of a to b is finite, and (ii) the last remainder (namely, the remainder which divides the immediately preceding remainder) is the Greatest Common Divisor of a, b.

*Elementary proof of Proposition X.2.*
We proceed with a proof by contradiction.

We assume that the lines a and b are commensurable, namely there is a line c and natural numbers m, n such that m > n and a = mc, b = nc.

By Propositions VII.1 & 2 the anthyphairesis of m to n is finite, say:

$m = k_0 n + e_1$ with $e_1 < n$,

$n = k_1 e_1 + e_2$ with $e_2 < e_1$,

… ,

$e_{p-1} = k_p e_p$.

But then

$mc = k_0 nc + e_1 c$ with $e_1 c < nc$,

$nc = k_1 e_1 c + e_2 c$ with $e_2 c < e_1 c$,

… ,

$e_{p-1} c = k_p e_p c$.

We set $d_k = e_k c$ for every $k = 0, 1, …, p$, and we have

$a = k_0 b + d_1$ with $d_1 < b$,

$b = k_1 d_1 + d_2$ with $d_2 < d_1$,

… ,

$d_{p-1} = k_p d_p$,

and thus the anthyphairesis of a to b is finite, a contradiction.

**5.2.4. The Theaetetean roots of the elementary proof of Proposition X.2**

The question that now arises is whether this elementary proof of Proposition X.2 is based on an idea of our own, or whether there is any indication that it might have been known to pre-Eudoxean mathematicians.

Indeed, the proof of Proposition VI.1 hinted by Aristotle, and employed for the proof of Theaetetus' analogue of Proposition VI.1 of the *Elements* (cf. Proposition 4.3.11), is based on the observation, rated as "evident at once" ("εὐθέως φανερὸν", Aristotle's *Topics*, 158b33), namely on the observation that:

if $Anth(a, b) = [k_0, k_1, k_2, …]$, then $Anth(ac, bc) = [k_0, k_1, k_2, …] = Anth(a, b)$:

ἔοικε δὲ καὶ ἐν τοῖς μαθήμασιν
ἔνια δι' ὁρισμοῦ ἔλλειψιν οὐ ῥαδίως γράφεσθαι, οἷον ὅτι
ἡ παρὰ τὴν πλευρὰν τέμνουσα τὸ ἐπίπεδον
ὁμοίως διαιρεῖ τήν τε γραμμὴν καὶ τὸ χωρίον.
τοῦ δὲ ὁρισμοῦ ῥηθέντος *εὐθέως φανερὸν* τὸ λεγόμενον·
τὴν γὰρ αὐτὴν ἀνταναίρεσιν ἔχει τὰ χωρία καὶ αἱ γραμμαί·
ἔστι δ' ὁρισμὸς τοῦ αὐτοῦ λόγου οὗτος. [Aristotle's *Topics* 158b29-35]

It also appears that in mathematics,

it is not easy to construct proofs of some things because of the deficiency of a definition, for instance that the line cutting a plane figure parallel to its side divides the area and the line similarly.

But once the definition has been stated, the proposition *is evident at once*.

For the areas and the lines have the same reciprocal subtraction [anthyphairesis]:

and this is the definition of 'same ratio'. [Translation by Smith, 1997]

But this is the exact observation, that

if Anth(m, n) = $[k_0, k_1, …, k_n]$, then Anth(a, b) = Anth(mc, nc) = Anth(m, n) = $[k_0, k_1, …, k_n]$,

on which our elementary proof of Proposition X.2 in Section 5.1.3, is based. Thus, certainly Theaetetus was in possession of the elementary proof of Proposition X.2.

Proposition X.2 is absolutely essential for the reconstruction of Theaetetus' incommensurabilities, including a correct reconstruction of Proposition X.9 of the *Elements*, as shown in Negrepontis (2018), Negrepontis, Farmaki and Kalisperi (2022), Negrepontis, Farmaki and Brokou (2024), and Negrepontis and Protopapas (2025, Sections 7.2.2 and 9). We emphasize that the proof of Theaetetus' anthyphairetic periodicity theorem makes use of the Pythagorean/Theodorus technique of Applications of Areas, and the novel element of the proof is essentially the use of the pigeon-hole argument.

## 5.3. Propositions X.10-16 (supplying useful properties), X.19-26 (concerning the definition, construction and properties of the medial line and area), X.29-30, X.84/85, X.85-90, X.47/48, X.48-53 and X.111-X.115 (concerning the definition and construction of the six apotome and the six binomial lines and the conjugation between the two hexads)

### 5.3.1. Proposition X.20

If a rational area is applied to a rational line, then it produces as breadth a line rational and commensurable in length with the line to which it is applied.

*Proof.* Let ab be the rational area applied to the rational line a. We now have to show that b is rational and commensurable in length with a. By Proposition I.46, we can describe on line a the square $a^2$, which is rational (Definition X.I.4). Since $a^2$, ab are rational areas, the Theaetetean analogue of VI.1 is valid (Proposition 4.3.11). Thus: $a^2/ab = a/b$. Since $a^2$ is commensurable with ab, a is commensurable with b (Proposition X.11). But a is rational, therefore b is also rational and commensurable in length with a (Definition X.I.1).

*Note.* There are no major differences in the proof Euclid demonstrates and our own in Proposition X.20. In fact, the only difference is the use of Proposition 4.3.11 from our Section

4 instead of VI.1, which is justifiable since the ratio involved is Theatetean as it consists of rational areas.

### 5.3.2. Definition (in Proposition X.21)

A line m is *medial* (with respect to a given line r) if there are lines a, b, such that each of a, b is rational, a, b are incommensurable to each other, and $m^2 = ab$.

### 5.3.3. Proposition X.22

If m is a medial line, a is a rational line, and $m^2 = ab$, then b is rational and incommensurable in length to a.

*Proof.*
Since m is medial line, by the definition given in X.21, there are lines $a_1$, $b_1$ rational, commensurable in square only, such that $m^2 = a_1 b_1$. Thus $ab = a_1 b_1$.

*Claim 1.* b is a rational line.
By definition of rational line, a, $a_1$ are commensurable in square, namely $Na^2 = Ma_1^2$ for some natural numbers M, N. Since $ab = a_1 b_1$, it follows from Proposition 4.2.3 that $Nb_1^2 = Mb^2$. Hence b is rational.

*Claim 2.* a, b are incommensurable in length.
Since $a_1$, $b_1$ are commensurable in square only, it follows that $a_1$, $b_1$ are incommensurable in length. Since $a_1$, $b_1$ are commensurable in square, Theaetetus' analogue of Proposition VI.1 (Proposition 4.3.11) is valid, thus $a_1/b_1 = a_1 b_1/b_1^2$, hence $b_1^2$, $a_1 b_1$ are incommensurable. Since b, $b_1$ are rational lines, it follows that $b^2$, $b_1^2$ are commensurable. Since $ab = a_1 b_1$, it follows that ab, $a_1 b_1$ are commensurable. Hence $b^2$, ab are incommensurable. Since a, b are rational lines, Theaetetus' analogue of VI.1 is valid, (Proposition 4.3.11) thus we have $a/b = ab/b^2$, hence a, b are incommensurable in length.

*Note.* The proof of Proposition X.22 differs significantly from its original, Euclidean, form in our reconstruction. We have avoided the use of Propositions VI.14 and VI.22 altogether as well as the preceding Lemma X.21/22. In fact, only Proposition 4.3.11 – our analogue to Proposition VI.1 in Euclid's *Elements* – is necessary for the proof to be completed and it is applicable in both cases since the ratios involved are Theatetean.

### 5.3.4. Proposition X.23

If m is a medial line, and n is a line commensurable to m, then n is a medial line.

*Proof.* Let a be a rational line. By Proposition I.44, there is a line b such that $m^2 = ab$. By Proposition X.22, since m is a medial line, it follows that b is rational and incommensurable to

a. By Proposition I.44, there is a line c such that $n^2 = ac$. Since m, n are commensurable lines, it follows that $m^2$, $n^2$ are commensurable areas, hence ab, ac are commensurable areas. Thus Theaetetus' analogue of Proposition VI.1 (Proposition 4.3.11) is valid, and thus $ab/ac = b/c$, hence b, c are commensurable. Since b is rational, it follows that c is rational, and since a, b are incommensurable, it follows that a, c are incommensurable. Thus $n^2 = ac$, and a, c are rational lines, commensurable in square only. By the definition in X.21, n is a medial line.

*Note.* There are no major differences in this proof from the one presented in the *Elements*. We stress again the fact that Proposition 4.3.11, the analogue of VI.1 from our reconstruction, is applicable here since the areas are commensurable and therefore generate a Theatetean ratio.

### 5.3.5. Definition (through their use in Proposition X.24)

A *medial* square $m^2$ is one whose side m is a *medial* line, and a *medial* rectangle ab is one equal to a square $m^2$ whose side m is a *medial* line.

### 5.3.6. Proposition X.26

A medial area does not exceed a medial area by a rational area.

*Proof.* Let the medial area A exceed the medial area B by a rational area R, thus $A - B = R$. Let r be a rational line set out. By Proposition I.44, there are lines a, b, such that $A = a \cdot r$, $B = b \cdot r$. Then the rectangles ar and br are medial as well. By Proposition X.22, since the line r is rational, the lines a and b are rational and incommensurable in length with r. Also $R = A - B = (a - b)r$. Since R is a rational area, so is the rectangle $(a - b)r$. By Proposition X.20 and since line r is rational, line $c = a - b$ is rational as well and commensurable in length with r. Therefore, b and c are rational lines incommensurable in length (Proposition X.13). Since b, c are rational lines Theaetetus' analogue of Proposition VI.1 (Proposition 4.3.11) is valid and thus, $b/c = b^2/bc$, hence $b^2$ is incommensurable with bc (Proposition X.11). Since $b^2$, $c^2$ are rational squares, the sum $b^2 + c^2$ is rational and commensurable with $b^2$, and since bc is commensurable with 2bc (Proposition X.6), the sum of the squares $b^2 + c^2$ is incommensurable with 2bc (Proposition X.13). Therefore, $b^2 + c^2 + 2bc$ is incommensurable with $b^2 + c^2$ (Proposition X.16). But, by Proposition II.4, $b^2 + c^2 + 2bc = (b + c)^2 = a^2$. Then $a^2$ is incommensurable with the rational sum $b^2 + c^2$, hence $a^2$ is irrational which means that a is irrational. But a was also proved to be rational, which is a contradiction.

*Note.* Once more, we must point out that in the use of Proposition 4.3.11, the lines involved are commensurable in square and hence the ratio produced is Theatetean.

### 5.3.7. Propositions X.29-30

Propositions X.29-30, although proven in Euclid's work, are essentially repeated in the constructions found in the proofs of Propositions X.48-53 and X.85-90, as noted in Section 5.3.9 below. We therefore consider them redundant and we will elaborate on their role as we encounter them further along.

**5.3.8. Definition of the apotome and binomial lines**

*X.36 Definition of the binomial line*
*X.47/48 Definition and Classification of the six binomial lines*
*X.73 Definition of the apotome line*
*X.84/85 Definition and Classification of the six apotome lines in terms of ζ, η, r, and of ζ, θ*

*Definition in Proposition X.73:*

A line γ is an apotome (with respect to a given line r) if there are lines ζ, η, such that
ζ > η, γ = ζ – η,
ζ, η are incommensurable,
$ζ^2, η^2, r^2$ are commensurable.

Two criteria are then used for the classification of the apotome lines:
(i) r compared to ζ, η [inspired from the anthyphairesis of a and b, when $a^2 = Nb^2$]
(ii) ζ compared to θ [inspired from the use of the proof of the palindromic periodicity theorem to prove the Pell property]

**5.3.9. Construction of the six binomial lines X.48-53 and the six apotome lines X.85-90**

*X.85-90. Construction of each of the six apotome lines*

The only place in Book X where apotome lines are expressly *constructed* is in the hexad of Propositions X.85-90. We notice that the apotome lines constructed there are somewhat less general than what one might expect from the definition of an apotome line. In order to explain this difference we introduce the following definition.

*Definition.* A line γ is a *simple apotome* if there are a non-square number N, lines a and b such that $a^2 = Nb^2$, and natural numbers p, q, λ, such that either λγ = pa – qb or λγ = qb – pa.
A simple apotome is *integral* if λ = 1, *fractional* if λ > 1.
The canonical choice for the given line of a simple apotome is the line b, corresponding to the one foot line of *Theaetetus* 147d5, or the line a, or any other line r commensurable in square to b, namely such that $a^2, r^2$ are commensurable areas.

In the *Elements*, the definition of a simple apotome is *not* explicitly given, but it is clear that Book X concerns simple apotomes, since all six kinds of apotome lines constructed in Propositions X.85-90, the only construction of apotome lines in Book X, are simple apotome lines.

Indeed, *for the construction of the line in Propositions X.85 – 87*, we let N be a non-square number, such that $N = q^2 – p^2$ (a method for finding such numbers is described in Lemma 1 in X.28/29), lines a, b such that $a^2 = Nb^2$, and construct the apotome line qb – a (with respect to given line r, take b in X.85, a in X.86, and incommensurable to both a and b in X.87).

*For the construction of the line in Propositions X.88 - 90*, we let N be a non-square number, such that $N = q^2 + p^2$ (a method for finding such numbers is described in Lemma 2 in X.28/29), lines a, b such that $a^2 = Nb^2$, and construct the apotome line a – qb (with respect to given line r, take a in X.88, b in X.89, and incommensurable to both a and b in X.90).
(Cf. Negrepontis, Farmaki, Brokou, 2024, Section 4.1.2).

*Note*. It must be emphasized that the general definition of an apotome line (in X.73) allows for lines of the form $\zeta – \eta$, where $\zeta^2 = Nb^2$, $\eta^2 = Mb^2$, for some two (different) non-square numbers N, M. Such an apotome line indeed appears in Proposition XIII.17 of the *Elements*: the side s of the face of the regular dodecahedron is the line $(1/3)\cdot(a – a')$, where $a^2 = 15d^2$, $(a')^2 = 3d^2$, d the diameter of the sphere in which the dodecahedron is comprehended.
Thus, Theaetetus in Book X is interested in the general apotome lines, so that a line, such as the one arising in the study of the dodecahedron, will be included in the definition of apotome lines and this appears to be the sole reason for adopting the general definition. His central interest however, as attested by all his explicit constructions in Book X, appears to be the *simple* apotome lines arising from the lines a, b, where $a^2 = Nb^2$, for some non-square number N.
At this point, we must also draw attention to the fact that only the Lemmas from Propositions X.29-30 are necessary for the proofs of Propositions X.85-90 (and X.48-53). However, the two Propositions find use later on, in the final part of Book X.

### 5.3.10. The conjugation between the apotome and binomial lines in Propositions X.112 - 114

*Proposition X.113*

The square on a rational line, if applied to an apotome, produces as breadth the binomial line the terms of which are commensurable with the terms of the apotome and in the same ratio; and further the binomial so arising has the same order as the apotome.

*Proof*. Let r be a rational line and $\gamma = \zeta – \eta$ an apotome with respect to the given line r. By Proposition I.44, construct a line $\delta$ such that the rectangle $\gamma\delta$ = the square $r^2$.

*Claim*. There are lines $\psi, \omega$ such that $\psi/\zeta = \omega/\eta$ is a commensurable ratio, $\delta = \psi + \omega$ is a binomial line with respect to the given line r, and and $\delta$ has the same order as $\gamma$.

*Proof of claim.*

By definition (in Proposition X.73), ζ and η are rational lines commensurable in square only. By Proposition I.44, construct a line φ such that the rectangle ζφ = the square $r^2$. But the square $r^2$ is rational, therefore the rectangle ζφ is also rational. Hence, by Proposition X.20, φ is rational and commensurable in length with ζ. Since rectangle ζφ = rectangle γδ, and ζ > γ, then δ > φ. It follows that rectangle ζ(δ – φ) = rectangle δ(ζ – γ). Since ζ/η is a Theaetetean ratio (Definition 4.3.1), and ζ – γ = η, we have, by Proposition 4.3.4 (the Theaetetean analogue of Proposition VI.16), ζ/η = δ/(δ – φ). Thus, δ/(δ – φ) is a Theaetetean ratio.

We next construct lines ω, χ, such that δ – φ = χ + ω and ζ/η = (φ + χ)/ω.

[Indeed, by Proposition VI.10 for the Theaetetean ratio ζ/η, construct lines ψ, ω, such that ψ + ω = δ, and ζ/η = ψ/ω. We set ψ – φ = χ].

Then, ψ = φ + χ, and δ = φ + χ + ω = ψ + ω. Since ζ and η are commensurable in square only and ζ/η = ψ/ω, it follows that ψ and ω are commensurable in square only. Since δ/(δ – φ) = ψ/ω, and δ – φ – ω = χ, δ – ψ = ω, by Proposition 4.3.9.b, the Theaetetean analogue for Proposition V.19, we have δ/(δ – φ) = ω/χ = (δ – ψ)/(δ – φ – ω), therefore ψ/ω = ω/χ; hence, since (informally) ψ/χ = (ψ/ω)·(ω/χ) = (ψ/ω)·(ψ/ω) = $ψ^2/ω^2$, we have ψ/χ = $ψ^2/ω^2$. Since ψ and ω are commensurable in square, then $ψ^2$ is commensurable with $ω^2$, hence ψ is commensurable in length with χ. Since ψ = φ + χ, it follows immediately that ψ is commensurable in length with φ. But φ is rational and commensurable in length with ζ, therefore ψ is also rational and commensurable in length with ζ. Since ζ/η = ψ/ω, by the *Alternando* property for Theaetetean ratios (Proposition 4.3.6), ζ/ψ = η/ω. But ζ is commensurable with ψ, therefore η is commensurable in length with ω. But ζ and η are rational lines commensurable in square only, therefore ψ and ω are also rational lines commensurable in square only. Therefore δ = ψ + ω is binomial.

In fact, by Propositions X.14, 48-53, it is easy to see that δ is a binomial line, the terms of which, ψ and ω, are commensurable with the terms ζ and η, respectively, of the apotome γ and in the same ratio, and δ has the same order as γ.

*Note.* The proof in this proposition is the first that makes use of propositions other than 4.3.11 from our reconstruction. Nevertheless, in all cases the ratios involved are Theatetean and, as such, permit the use of propositions from Theaetetus' theory of ratios, as it is presented in Section 4. In addition, the use of Proposition VI.10 is in this case permissible, even though the proof of VI.10 makes use of VI.2 and Eudoxus' condition is involved in the proof of VI.2. However, only the converse part of the statement in Proposition VI.2 requires the use of Eudoxus' condition and we do not need the converse here.

## 5.4. Definition and basic properties of the 13 alogos lines

### 5.4.1. Definitions

*X.21 definition of medial line* (*already given in 5.3.2*)
*X.36-41 Definition of the six kinds of additive alogos lines and proof of their alogia*
*X.73-78 Definition of* the *six kinds of subtractive alogos lines and proof of their alogia*

*Definition in Proposition X.73:*

Recall that the definition of an apotome has already been given in Section 5.3.8.

*Definition in Proposition X.75:*

A line Ω is a second apotome of a medial line if there are lines Φ, Ψ, such that
Φ > Ψ, Ω = Φ – Ψ,
Φ is a medial line, Ψ is a medial line,
Φ, Ψ commensurable in square only, and
ΦΨ is a medial area.

### 5.4.2. Alogia of the alogoi lines

*Propositions X.36-41 for the construction of the six additive alogos lines*
*Propositions X.73-78 for the construction of the six subtractive alogos lines*

The proofs of Propositions X.73, the alogia of *the apotome line*, and X.75, the alogia of *the second apotome of a medial line*, will be given in detail. All proofs in these two hexads are similar.

#### 5.4.2.1. Proposition X.73

Every apotome line is an alogos line.

*Proof.* Let γ = ζ – η be an apotome line wrt a given line r, namely ζ > η, ζ, η are rational lines, and ζ, η are commensurable in power/square only. Then ζ, η are incommensurable in length. Since ζ, η are commensurable in square, the Theaetetus' analogue of VI.1 (Proposition 4.3.11) is valid, thus ζ/η = ζ$^2$/ηζ, hence ζ$^2$, ηζ are incommensurable. On the other hand, ζ$^2$, η$^2$, hence ζ$^2$ + η$^2$ as well, are commensurable and ηζ, 2ηζ are commensurable. By Proposition II.7, γ$^2$ = (ζ – η)$^2$ = ζ$^2$ + η$^2$ – 2·ζ·η, hence ζ$^2$ + η$^2$, γ$^2$ are incommensurable, hence r$^2$, γ$^2$ are incommensurable, hence the apotome line γ is an *alogos* line.

#### 5.4.2.2. Proposition X.75

Every second apotome of a medial line is an alogos line.

*Proof.* Let Ω be a second apotome of a medial line, Ω = Φ – Ψ and Φ, Ψ are incommensurable. Since Φ, Ψ are commensurable in square, Theaetetus' version of Proposition VI.1 (Proposition 4.3.11) is valid, thus $\Phi/\Psi = \Phi^2/\Phi\Psi$, hence $\Phi^2$, $\Phi\Psi$ are incommensurable; but $\Phi^2$, $\Psi^2$, hence $\Phi^2 + \Psi^2$ as well, are commensurable and $\Phi\Psi$, $2\Phi\Psi$ are commensurable. Hence $\Phi^2 + \Psi^2$, $2\Phi\Psi$ are incommensurable. By Proposition I.44, construct line $\zeta$, such that $\Phi^2 + \Psi^2 = r\cdot\zeta$. By Proposition I.44, construct line $\eta$, such that $2\cdot\Phi\cdot\Psi = r\cdot\eta$. Set $\gamma = \zeta - \eta$. By Proposition II.7, $\Omega^2 = (\Phi - \Psi)^2 = \Phi^2 + \Psi^2 - 2\Phi\Psi = r\cdot\zeta - r\cdot\eta = r\cdot(\zeta - \eta) = r\,\gamma$. Since $\Phi$, $\Psi$ are medial lines, it follows that $\Phi^2$, $\Psi^2$ are medial areas, hence, by Proposition X.15, $\Phi^2 + \Psi^2$ is a medial area. By Proposition X.22, $\zeta$ is a rational line, and r, $\zeta$ are commensurable in square only. Since $\Phi\Psi$ is a medial area, then $2\Phi\Psi$ is a medial area. By Proposition X.23, $\eta$ is a rational line, and r, $\eta$ are commensurable in square only. Since $\zeta$, $\eta$ are rational lines, Theaetetus' version of Proposition VI.1 (Proposition 4.3.11) is valid and we have $\zeta/\eta = r\zeta/r\eta = (\Phi^2 + \Psi^2)/2\Phi\Psi$. Thus, $\zeta$, $\eta$ are rational and incommensurable to each other, hence line $\gamma$ is an apotome line, and by Proposition X.73, proved above, line $\gamma$ is an alogos line. Since $\Omega^2 = r\gamma$ and r is a rational line, it follows that $\Omega$ is an alogos line.

*Note.* Our proofs of Propositions X.73 and X.75 are similar to Euclid's. The application of Proposition 4.3.11 is permissible wherever we use it in both proofs, since the magnitudes involved generate Theatetean ratios.

**5.4.3. Uniqueness of the alogoi lines**

*Propositions X.79-84 for the six subtractive alogos lines,*
*Propositions X.42-47 for the six additive alogos lines*

The proof of Proposition X.82, the uniqueness of the minor alogos line, will be given in detail. All proofs in these two hexads are similar. The proofs of the twelve uniqueness propositions are based on Propositions X.22 and X.26.

*Proposition X.82.* To a minor line only one line can be annexed which is incommensurable in square with the whole and which makes, with the whole, the sum of the squares on them rational but twice the rectangle contained by them medial.

*Proof.* Let $\Omega$ be the minor line, and let $\Psi$ be an annex to $\Omega$. Then $\Phi = \Omega + \Psi$ and $\Psi$ are lines incommensurable in square which make the sum of the squares on them rational, but twice the rectangle contained by them medial (definition in Proposition X.76). By Proposition II.7, $\Omega^2 = (\Phi - \Psi)^2 = \Phi^2 + \Psi^2 - 2\Phi\Psi$. Now, let Y be another line which can be annexed to $\Omega$ and consider the case Y > $\Psi$. Then $\Phi = \Omega + Y$ and Y are lines incommensurable in square which make the sum of the squares on them rational, but twice the rectangle contained by them medial. By Proposition II.7, $\Omega^2 = (\Phi - Y)^2 = \Phi^2 + Y^2 - 2\Phi Y$. Therefore $\Phi^2 + Y^2 - 2\Phi Y = \Phi^2 + \Psi^2 - 2\Phi\Psi$, hence $2\Phi Y - 2\Phi\Psi = \Phi^2 + Y^2 - (\Phi^2 + \Psi^2)$ which is a contradiction to Proposition X.26, since by the last equality we have proved that the medial area $2\Phi Y$ exceeds the medial area $2\Phi\Psi$ by the rational area $\Phi^2 + Y^2 - (\Phi^2 + \Psi^2)$.

*Note.* Our proof follows Euclid's reasoning in this proposition. No ratios are involved in this proposition's proof.

### 5.4.4. Invariance under commensurability of the alogoi lines

*Propositions X.66 - 70 for the additive alogoi lines,*
*Propositions X.103 - 107 for the subtractive alogoi lines.*

The proof of Proposition X.105, the invariance under commensurability of the alogos line minor will be given in detail. All proofs in these two hexads are similar.

*Proposition X.105.*
If $\Omega$ is a minor line and $\omega$ is a line commensurable with $\Omega$, then $\omega$ is a minor line.

*Proof.* Since $\Omega$ is a minor line, then $\Omega = \Phi - \Psi$, such that $\Phi, \Psi$ are incommensurable in square, $\Phi^2 + \Psi^2$ is a rational area, and $\Phi\Psi$ is a medial area.
Let $\omega$ be commensurable to $\Omega$, say $\omega/\Omega = m/n$ for some natural numbers m, n. We define $\psi, \varphi$ by the conditions $\omega/\Omega = \psi/\Psi$, $\varphi = \omega + \psi$. Thus, $\omega = \varphi - \psi$, and $\varphi/\Phi = \omega/\Omega = \psi/\Psi = m/n$ (Proposition 4.3.9.a).
Since $\Phi, \Psi$ are incommensurable in square, $\varphi, \Phi$ are commensurable, and $\psi, \Psi$ are commensurable, it follows that $\varphi, \psi$ are incommensurable in square.
Since $\varphi/\Phi = \psi/\Psi = m/n$, it follows that $\varphi^2/\Phi^2 = \psi^2/\Psi^2 = m^2/n^2$, hence $(\varphi^2 + \psi^2)/(\Phi^2 + \Psi^2) = m^2/n^2$ (Proposition 4.3.9.a for areas).
Since $\Phi^2 + \Psi^2$ is a rational area and $\varphi^2 + \psi^2$ is an area commensurable to it, it follows that $\varphi^2 + \psi^2$ is a rational area.
Since $\varphi/\Phi = \psi/\Psi = m/n$, it follows that $\varphi\psi/\Phi\Psi = m^2/n^2$. Since $\Phi\Psi$ is a medial area and $\varphi\psi$ is an area commensurable to it, it follows from Proposition X.23 (and its porism) that $\varphi\psi$ is a medial area. It follows that $\omega$ is a minor line.

*Note.* Our proof is much simpler than Euclid's in this proposition and avoids the use of many propositions on ratios. Once again though, all ratios are Theatetean and, as such, warrant the use of propositions from Theatetus' theory of ratios.

### 5.5. From apotome/binomial lines to alogoi lines (X.91-96, X.54-59) and conversely (X.97-102, X.60-65), *with auxiliary Propositions X.27-28, 31-35 and the help* of Application of Areas in Defect

The proofs of the propositions in this subsection, in which the Pythagorean Application of Areas in defect is employed, make use of substantially new arguments in order to be derivable from the restricted Theaetetean theory of proportion. Thus, Lemma X.32/33 receives an elementary proof with no theory of proportion.

### 5.5.1. Lemma X.32/33

Let ABC be a right-angled triangle having the angle A right, and let the perpendicular AD be drawn. Then the rectangle CB by BD equals the square on BA, the rectangle BC by CD equals the square on CA, the rectangle BD by DC equals the square on AD, and the rectangle BC by AD equals the rectangle BA by AC.

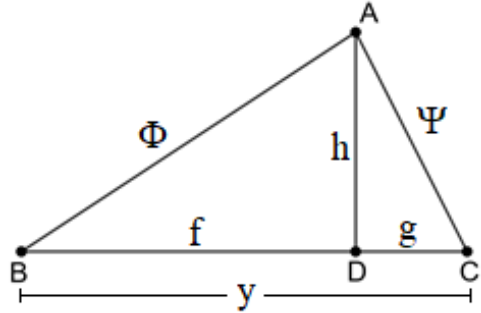

*Figure 1: right-angled triangle ABC, as described in the proof of Lemma X.32/33*

The elementary *proof* does not need Eudoxus' theory of ratios of magnitudes which Euclid uses, but follows from Proposition II.4 and the Pythagorean theorem.

*Proof.* Let $AB = \Phi$, $AC = \Psi$, $BC = y$, $AD = h$, $BD = f$, $DC = g$. We must prove that $\Phi^2 = yf$, $\Psi^2 = yg$, $h^2 = fg$ and $\Phi\Psi = yh$.
By the Pythagorean Theorem and II.4, $\Phi^2 + \Psi^2 = y^2 = (f + g)^2 = f^2 + g^2 + 2fg$. Also, $\Phi^2 = h^2 + f^2$ and $\Psi^2 = h^2 + g^2$. Then $f^2 + g^2 + 2h^2 = f^2 + g^2 + 2fg$, hence $h^2 = fg$.
Then, $\Phi^2 = h^2 + f^2 = fg + f^2 = (f + g)f = yf$, thus $\Phi^2 = yf$, and $\Psi^2 = h^2 + g^2 = fg + g^2 = (f + g)g = yg$, thus $\Psi^2 = yg$.
Finally, $\Phi\Psi = hy$.

### 5.5.2. The use of Application of Areas in defect in Book X of the *Elements*

All the propositions in this Section employ the Application of Areas in defect. The only proposition that Euclid has proved and to which he can refer for the proof of these propositions is Proposition VI.28, a proposition that depends heavily on Eudoxus' theory of ratios and on Eudoxus' condition. Thus, the use of Proposition VI.28, as Heath suggests in the proofs of Propositions X.17 and X.18, is problematic for our reconstruction. However, the only proposition needed throughout Book X is the Pythagorean Application of Areas in defect, which would undoubtedly have been familiar to Theaetetus. We therefore find a solution to the problem by interpolating between Propositions II.5 and II.6 of Book II what we will denote by Proposition II.5/6, namely a Pythagorean Application of Areas in defect.

Indeed, the existence of such a proposition has long had its advocates, of which we mention here Heath himself (1926) and van der Waerden (1954). Furthermore, in a complete and thorough study of Book II of the *Elements* leading to the restoration of the Book to its original, Pythagorean form Negrepontis, Farmaki (2025) have convincingly argued that additional propositions of Book II must have been removed by Euclid since they were in essence repeated in other Books of the *Elements*. Therefore, in our reconstruction of Book X we may replace proposition VI.28 with:

*Proposition II.5/6. (Pythagorean Application of Areas in Defect)*

For a and m lines, with $a/2 > m$, to construct a line x, such that $x(a - x) = m^2$.

### 5.5.3. Propositions X.17 and X.18

*Proposition X.17*

If there be two unequal lines, and to the greater there be applied a parallelogram equal to the fourth part of the square on the less and deficient by a square figure, and if it divide it into parts which are commensurable in length, then the square on the greater will be greater than the square on the less by the square on a line commensurable with the greater. And, if the square on the greater be greater than the square on the less by the square on a line commensurable with the greater, and there be applied to the greater a parallelogram equal to the fourth part of the square on the less and deficient by a square figure, it will divide it into parts which are commensurable in length.

*Proof.*

Let r be a rational line and let $\zeta, \eta$ be two unequal lines with $\zeta > \eta$.

*Part 1*: [*Construction process, depicted in figure 2*]

Bisect η and apply to ζ a parallelogram equal to the square on η/2 and deficient by a square figure, and let it be the rectangle ζ− x by x (Proposition II.5/6). Describe the circle on AB = ζ. Draw ED at right angles to AB. Join AD = Φ and DB = Ψ. Then D is a right angle, by Proposition III.31. By hypothesis: $(ζ − x)x = η^2/4$. By Lemma X.32/33, $(ζ − x)x = DE^2$ and therefore, DE = η/2. Also, CE = (ζ/2) – x = (ζ – 2x)/2. By the Pythagorean Theorem, $ζ^2/4 = η^2/4 + CE^2$. Then, $ζ^2/4 = η^2/4 + (ζ − 2x)^2/4$. And by quadrupling we have that: $ζ^2 = η^2 + (ζ − 2x)^2$. Let ζ − 2x = θ ⇔ θ/2 = ζ/2 − x = CE.

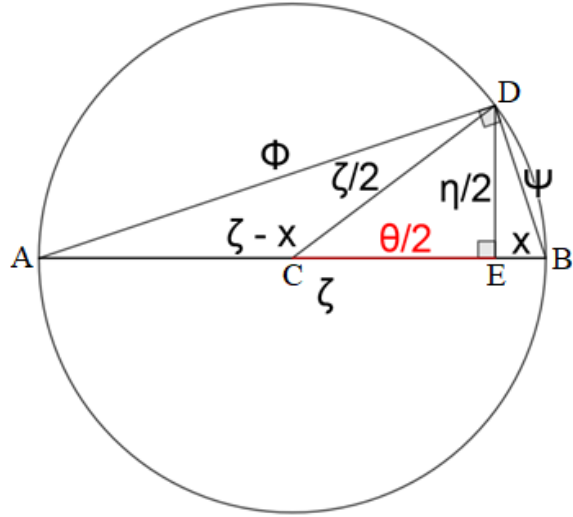

*Figure 2: construction of lines θ, Φ and Ψ from lines ζ and η.*

*Part 2*: [*proof of same type of commensurability between ζ – x, x and ζ, θ*]

We now have Figure 2 from the construction explained in Part 1. Let x and ζ − x be commensurable. Then, m(ζ – x) = nx, for some natural numbers m, n and thus, mζ = (m + n)x, and we have that m < n. But then x and ζ are commensurable (Proposition X.15). Also, x and 2x are obviously commensurable (Proposition X.6). Therefore, ζ and 2x are commensurable (Proposition X.12), since 2mζ = 2(m + n)x ⇔ (2m)ζ = (m + n)2x, and thus, ζ and θ = ζ − 2x are commensurable (Proposition X.15).
Indeed, (m + n)ζ – (2m)ζ = (m + n)ζ – (m + n)2x ⇔ (n – m)ζ = (m + n)(ζ – 2x) ⇔ (n – m)ζ = (m + n)θ

*Conversely*:

We know that ζ and ζ − 2x are commensurable. Then we can find natural numbers m, n with m < n such that mζ = n(ζ − 2x). Therefore, by Proposition X.15, (n – m)ζ = n(2x) and hence ζ and 2x are commensurable. This, in turn, means that ζ and x are commensurable, that is, (n – m)ζ = (2n)x and, finally, x and ζ − x are commensurable (by Proposition X.15) since (n – m)ζ – (n – m)x = (2n)x – (n – m)x ⇔ (n – m)(ζ – x) = (m + n)x.

*Proposition X.18*

If there be two unequal lines, and to the greater there be applied a parallelogram equal to the fourth part of the square on the less and deficient by a square figure, and if it divide it into parts which are incommensurable, the square on the greater will be greater than the square on the less by the square on a line incommensurable with the greater. And, if the square on the greater be greater than the square on the less by the square on a line incommensurable with the greater, and if there be applied to the greater a parallelogram equal to the fourth part of the square on the less and deficient by a square figure, it divides it into parts which are incommensurable.

*Proof.* Let r be a rational line. Let ζ, η be two unequal lines with ζ > η. Consider Figure 2 from Proposition X.17 but let x and ζ − x be incommensurable this time. We will show that ζ and ζ − 2x are also incommensurable.

Indeed, if ζ and ζ − 2x were commensurable, then x and ζ – x would be commensurable as well (Proposition X.17, converse part) and we would have a contradiction.

Therefore ζ and ζ − 2x = θ are incommensurable.

*Conversely:*

Let now ζ and ζ − 2x be incommensurable. We will show that x and ζ − x are commensurable. Indeed, by working in a way identical to what we have just demonstrated, if x and ζ − x were commensurable, then ζ and ζ – 2x would be commensurable as well (Proposition X.17) and we would have a contradiction.

Therefore x and ζ − x are incommensurable.

*Note.* We follow Euclid's path in the proof of X.17, but the proof of X.18 is quite different in our reconstruction and makes use only of its preceding proposition. Heath indicates that Proposition VI.28 is necessary in Euclid's proof, yet no ratios need to be involved in any of the two proofs. In fact, only a use of Proposition II.5/6, a Pythagorean application of areas, is necessary for the proof to be completed.

**5.5.4. Propositions X.27-28 and X.31-35**

The seven Propositions X.27-28 and X.31-35 are useful in proving the existence and constructability of the alogoi lines. The process followed in the proof of some of them is also quite insightful in providing a means for the proof of the most difficult propositions to follow, such as X.94 and X.100.

*Proposition X.27*
To find medial lines commensurable in square only which contain a rational rectangle.

*Proof.* Let ζ, η be rational lines, commensurable in square only. Then the rectangle ζη is medial. (Proposition X.21). Place ζ, η on a line and describe the semicircle with diameter ζ + η. Draw segment h, inscribed within the semicircle, perpendicular to ζ + η at their common endpoint (Proposition I.11). By joining the not common endpoints of segments ζ, η and h a triangle is formed, and since it is inscribed in the semicircle, it is a right-angled triangle (Proposition III.31). By Lemma X.32/33, $h^2 = ζη$ and then $h^2$ is a medial square (Proposition X.23, corollary). Then the line h is medial. Let it be contrived that ζ is to η as h is to x. Since ζ and η are lines commensurable in square only, h and x are lines commensurable in square only as well (Proposition X.11). Since h is a medial line, x is a medial line also (Proposition X.23).

Then by the Theatetean variant of Proposition VI.16 (Proposition 4.3.4, applicable here since both ratios are of lines commensurable in square) $\zeta/\eta = h/x \Leftrightarrow \zeta x = \eta h$. Multiplying both sides by h, we have: $\zeta xh = \eta h^2$. And since $h^2 = \zeta \eta$, we then arrive at: $\zeta xh = \eta \zeta \eta$. Simplifying ζ from both sides, we receive: $xh = \eta^2$. Therefore the rectangle xh is rational, since the square $\eta^2$ is rational.

Thus, we have found x and h, two medial lines commensurable in square only which contain a rational rectangle.

*Note:* In the original Euclidean proof of Proposition X.27, there is a use of the *Alternando* property for ratios which results in a ratio that is not Theatetean. However, this is circumvented in our proof since the mean proportional is found through the use of Lemma X.32/33 and thus without use of any propositions from Book VI.

### 5.5.5. From the six apotome and the six binomial lines to the six additive alogos lines X.54-59, and to the six subtractive alogos lines X.91-96

*Proposition X. 94*
If an area be contained by a rational line and a fourth apotome, the "side" of the area is minor.

*Proof.* Let r be a rational line and γ a fourth apotome. Then $\gamma = \zeta - \eta$, where ζ, r are commensurable in length and ζ, η are commensurable in square only. Therefore r, η are commensurable in square only. Also $\zeta^2 = \eta^2 + \theta^2$ with ζ, θ incommensurable in length. We need to show that the line ω, such that $\omega^2 = r\gamma$, is minor. Recall Figure 2 from Proposition X.17. By Lemma X.32/33, $\Phi\Psi = \zeta\eta/2$, which means that ΦΨ is medial. $\Phi^2 + \Psi^2 = \zeta^2$ (by the Pythagorean Theorem), which means that the sum of $\Phi^2$ and $\Psi^2$ is rational. By Proposition X.18, $\zeta - x$ and x are incommensurable, since ζ and θ are incommensurable. Then rectangles $\zeta(\zeta - x)$ and $\zeta x$ are incommensurable, which means that the squares $\Phi^2$ and $\Psi^2$ are incommensurable, so Φ and Ψ are incommensurable in square. We define then the minor $\Omega = \Phi - \Psi$.
$\Omega^2 = (\Phi - \Psi)^2 = \Phi^2 + \Psi^2 - 2\Phi\Psi = \zeta(\zeta - x) + \zeta x - 2\zeta\eta/2 = \zeta^2 - \zeta\eta = \zeta(\zeta - \eta) = \zeta\gamma$.
But ζ, r are commensurable in length and this means that $M\zeta = Nr$ for natural numbers M, N. Then $r/\zeta = M/N$. But $M\Omega^2 = M\zeta\gamma = Nr\gamma = N\omega^2$, which in turn means that $\omega^2/\Omega^2 = M/N$, and thus $\Omega^2$ and $\omega^2$ are commensurable in square. We define φ, ψ such that: $\varphi^2/\Phi^2 = M/N$ and $\psi^2/\Psi^2 = M/N$. Then we have: $\varphi^2/\Phi^2 = M/N = r/\zeta$. By Proposition VI.1 (Proposition 4.3.11 of our reconstruction) we have $r/\zeta = r(\zeta - x)/\zeta(\zeta - x) = r(\zeta - x)/\Phi^2$. Then: $\varphi^2/\Phi^2 = r(\zeta - x)/\Phi^2$ and therefore, $\varphi^2 = r(\zeta - x)$ (Proposition V.9 or 4.3.12 of our reconstruction). Similarly, we can show that $\psi^2 = rx$.

*Lemma.* $\varphi\psi = r(\eta/2)$.

*Proof of Lemma.* By Proposition I.44 we can find η' such that φψ = r(η'/2). Then (ζ – x)r(η'/2) = φ²(η'/2), and also (ζ – x)r(η'/2) = (ζ – x)φψ. Thus φ²(η'/2) = (ζ – x)φψ, and we therefore have φ(η'/2) = (ζ – x)ψ. Furthermore φψ(η'/2) = [φ(η'/2)]ψ = (ζ – x)ψ² = (ζ – x)rx, and also φψ(η'/2) = r(η'/2)(η'/2) = r(η'/2)². Thus, (ζ – x)rx = r(η'/2)², therefore (ζ – x)x = (η'/2)². But (ζ – x)·x = (η/2)². Therefore (η'/2)² = (η/2)², which leads us to η' = η. But then φψ = r(η/2).

Since φ²/Φ² = M/N = ψ²/Ψ² we have (φ² + ψ²)/(Φ² + Ψ²) = M/N (analogue to Proposition V.12 or of 4.3.9.a for areas). Then the area φ² + ψ² is commensurable with Φ² + Ψ², so both are rational. Since φψ/ΦΨ = r/ζ = M/N, the area φψ is commensurable with ΦΨ, therefore so are the areas 2φψ and 2ΦΨ. Since ΦΨ is medial, so is 2φψ. Moreover, since Φ, Ψ are incommensurable in square, φ, Φ are commensurable, and ψ, Ψ are commensurable, it follows that φ, ψ are incommensurable in square. We then define ω = φ – ψ, therefore ω is a minor. Also ω² = (φ – ψ)² = φ² + ψ² + 2φψ = r(ζ – x) + rx – rη = rζ – rη = r(ζ – η) = rγ. Therefore, ω is the line we are looking for.

*Note.* Our proofs for the hexad of Propositions X.91-96 vary significantly from the corresponding proofs Euclid presents in his book. The preceding proof of Proposition X.94 is indicative of our differentiation. Moreover, in this proposition's proof we encounter several uses of propositions from Theaetetus' theory of ratios. Again though, all ratios we encounter are Theatetean and we are therefore justified in using proposals from Theaetetus' theory of ratios.

### 5.5.6. From alogoi lines back to binomial and apotome lines Propositions X.97-102, X.60-65

This hexad contains the proposition with the most delicate proof. It is the only hexad (of propositions on apotomes), whose proofs make use of a Theaetetean ratio that has eventually periodic anthyphairesis, not because it is a commensurable in square ratio, but because it satisfies an Application of Areas in defect.

*Proposition X.100*
If Ω is a minor line, r is a rational line and γ is the line (constructed by Proposition I.44) such that Ω² = γr, then the line γ is a fourth apotome (defined in X.84/85, cf. 5.3.9).

*Proof.* Let Ω = Φ – Ψ be a minor (defined in Proposition X.76) and r a rational line. Then the line γ such that Ω² = γr is a fourth apotome. Since Ω is a minor, lines Φ, Ψ are incommensurable in square, Φ² + Ψ² is rational, and 2ΦΨ is medial. By Proposition I.44 construct lines ζ, η, such that Φ² + Ψ² = ζr, and 2ΦΨ = ηr. Then ζ > η. Define line γ by γ = ζ – η. Then, since Φ² + Ψ² is rational, ζr is rational as well, and since r is rational, by Proposition X.20, line ζ is rational and commensurable in length with r; and since 2ΦΨ is medial, ηr is medial as well, and since r is rational, by Proposition X.22, line η is rational and commensurable in square only with r. Then ζ, η are rational lines commensurable in square only. Also: $\Omega^2 = (\Phi - \Psi)^2 = \Phi^2 + \Psi^2 - 2\Phi\Psi = \zeta r - \eta r = (\zeta - \eta)r = \gamma r$. Thus, γ is an apotome (defined in Proposition X.73).

By application of areas in defect (Proposition II.5/6), line x is constructed such that $x(\zeta - x) = \eta^2/4$. In order to prove that γ is a fourth apotome, we need to show that ζ and $\theta = \zeta - 2x$ are incommensurable in length and, by Proposition X.18, it will be sufficient to prove that $\zeta - x$ and x are incommensurable.

By Proposition I.44, there is a line ψ such that $\Psi^2 = \psi r$. Then $\Phi^2 = (\zeta - \psi)r$. By Lemma X.32/33, $h^2 = fg$. Then $h^2 y = fgy$. But, also by Lemma X.32/33, $hy = \Phi\Psi = (\eta/2)r$, and $gy = \Psi^2 = \psi r$. Hence $(\eta/2)rh = \psi rf$, hence $(\eta/2)h = \psi f$, hence $(\eta/2)hy = \psi fy$. But, once more by Lemma X.32/33, $hy = \Phi\Psi = (\eta/2)r$ and $fy = \Phi^2 = (\zeta - \psi)r$. Hence $(\eta^2/4)r = \psi(\zeta - \psi)r$, and, finally, $\eta^2/4 = \psi(\zeta - \psi)$.

It follows that $\psi = x$.
[Indeed $\zeta^2 > \eta^2$ and thus $\zeta^2 - \eta^2 = \zeta^2 - 4x(\zeta - x) = (\zeta - 2x)^2$ and similarly $\zeta^2 - \eta^2 = \zeta^2 - 4\psi(\zeta - \psi) = (\zeta - 2\psi)^2$. Therefore, $(\zeta - 2x)^2 = (\zeta - 2\psi)^2$ and this leads to $\psi = x$.]
Thus, x and $\zeta - x$ are incommensurable. It follows from Proposition X.18 that ζ and θ are incommensurable. Thus, finally, line γ is a *fourth apotome*.

*Note.* In our proof of X.100, Proposition II.5/6 replaces any need of Eudoxus' condition and the elementary Lemma X.32/33 replaces the illicit, in the Theaetetean theory of ratios of magnitudes, uses of Proposition VI.1 in the *Elements* in order to prove that $\psi(\zeta - \psi) = \eta^2/4$, namely the ratios $\Phi/\Psi = \Phi^2/\Phi\Psi = (\zeta - \psi)r/(\eta/2)r = (\zeta - \psi)/(\eta/2)$ and $\Phi/\Psi = \Phi\Psi/\Psi^2 = (\eta/2)r/\psi r = (\eta/2)/\psi$. We would therefore conclude that $(\zeta - \psi)/(\eta/2) = (\eta/2)/\psi$, hence, by VI.16, $(\zeta - \psi)\psi = \eta^2/4$. The problem in Euclid's proof is apparent in the ratio Φ/Ψ, where the lines Φ, Ψ are incommensurable in square. In our proof this is avoided, as this ratio is not Theaetetean. Therefore, for the reasons mentioned directly above, we could not have accepted Euclid's proof of Proposition X.100 as it was written in the *Elements*. Instead, we have presented an entirely different, novel argument to prove X.100 in a manner that is befitting to our reconstruction proposal and in which only Theaetetean ratios appear.

## Section 6. The structure of Book X and a restoration of the correct historical order of the Books of the *Elements*

### 6.1. The structure of the restored Book X of the *Elements*

Summarizing, we will present in Table 1 the essential structure and order of Book X of the *Elements* that results from our reconstruction and restoration in the present work, an order that does not coincide with the order in the *Elements*.

In our reconstruction, we placed Theaetetus' theory of proportion (all propositions described in Section 4 of the present paper) in the first part of Book X and hence these propositions occupy column [1] of our resulting Table 1, which depicts the structure of Book X.

The second part then of Book X (column [2] of Table 1) consists of Propositions X.1-9 (related to Theaetetus' first incommensurability theorem reported in the *Theaetetus* 147d-e and described philosophically in Plato's *Sophist*) where the proofs of Propositions X.2 and X.9 have been modified as demonstrated in Section 5.2.

The propositions we have placed in the next part of Book X, occupying column [3] of Table 1, consist of the set of Propositions X.10-16 and X.19-26, which provide us with two definitions (medial line and area), their construction process and some useful auxiliary results, and also of the set of Propositions X.29-30, X.47/48, X.48-53, X.84/85, X.85-90 and X.111-115. It is in this latter set that we find the definition (X.84/85, X.47/48) and construction process (X.29-30, X.85-90, X.48-53) of the apotome and binomial lines and the conjugation between them (X.111-115).

We subsequently construct the alogoi lines of Book X and we have therefore placed all corresponding propositions in column [4] of Table 1. Specifically, in the fourth part of Book X we have the propositions that give us the definition of each alogos line and the proof of its alogia (X.73-78 and X.36-41), the propositions regarding the uniqueness of each alogos line (X.79-84 and X.42-47) and the propositions establishing each alogos line's invariance under commensurability (X.103-107 and X.66-71).

In the final, fifth, part of Book X, which takes up column [5] of Table 1, we have placed the propositions that provide us with the link between each apotome and binomial line and its corresponding alogos line. It is in this part that we find Propositions X.91-96 and X.54-59, in which we have a process of constructing an alogos line from an apotome/binomial line (used for the regular icosahedron in Book XIII), while the converse process is demonstrated in Propositions X.97-102 and X.60-65 (related to the proof of the Pell property, with the help of Propositions X.27-28, X.31-32 that lead to the application of areas in defect X.33-35, and related to Plato's *Meno* 86-87). The reader will have noticed that the proofs of the propositions contained in the final part of Book X were the most demanding and delicate in their process and were thus the ones that required major alterations from their original, Euclidean form.

*Note.* In Table 1, below, we indicate within brackets […] the propositions we proved in detail in Section 5, as representatives of a corresponding hexad or group. The proofs of the propositions in the same hexad, whether of apotome or binomial lines, are entirely similar and in our proofs we always took an apotome line as a representative of each hexad.

Table 1: The structure of Book X of the *Elements*

| Order of propositions in Book X of the *Elements* | [1] Theaetetus' theory of ratios of magnitudes | [2] Theaetetus' theorem on anthyphairetic periodicity | [3] Medial, apotome, binomial lines; theorem on palindromic periodic anthyphairesis of surds | [4] Definition and properties of alogoi lines | [5] Relation of apotome/binomial and alogoi lines by means of anthyphairetic periodicity of Application of Areas in Defect, Book XIII, Pell |
|---|---|---|---|---|---|
| Propositions X.1-35 ||||||
|  | Section 2 of present paper |  |  |  |  |
| X.1-9 |  | X.1-9 [X.2, X.9] |  |  |  |
| X.10-16 |  |  | X.10-16 |  |  |
| X.17-18 |  |  |  |  | X.17-18 [X.17-18] |
| X.19-26 |  |  | X.19-26 [X.22-23, X.26] |  |  |
| X.27-32 |  |  | X.29-30 |  | X.27-28, X.31-32 [X.27] |
| X.33-35 |  |  |  |  | X.33-35 [X.32/33] |
| Propositions on additive (binomial and alogoi) lines, Propositions X.36-72 ||||||
| X.36-41 |  |  |  | X.36-41 |  |
| X.42-47 |  |  |  | X.42-47 |  |
| X.47/48 |  |  | X.47/48 |  |  |
| X.48-53 |  |  | X.48-53 |  |  |
| X.54-59 |  |  |  |  | X.54-59 |
| X.60-65 |  |  |  |  | X.60-65 |
| X.66-70, 71-72 |  |  |  | X.66-70, 71-72 |  |
| Propositions on subtractive (apotome and alogoi) lines, Propositions X.73-110 ||||||
| X.73-78 |  |  |  | X.73-78 [X.73, X.75] |  |
| X.79-84 |  |  |  | X.79-84 [X.82] |  |
| X.84/85 |  |  | X.84/85 |  |  |
| X.85-90 |  |  | X.85-90 [X.85-90] |  |  |
| X.91-96 |  |  |  |  | X.91-96 [X.94] |
| X.97-102 |  |  |  |  | X.97-102 [X.100] |
| X.103-107, 108-110 |  |  |  | X.103-107, 108-110 [X.105] |  |
| Propositions on the conjugation of apotome and binomial lines, Propositions X.111-115 ||||||
| X.111-115 |  |  | X.111-115 [X.113] |  |  |

## 6.2. Restoration of the historical order of the thirteen Books of the *Elements*

The fact that Euclid does not include any mention of Theaetetus' theory of proportion for magnitudes forces him to an early appearance of Eudoxus' theory of ratios for magnitudes, and thus to a distortion of the true historical sequence. Once this is realized, there is no serious problem to proceed with the restoration of the historical order of the thirteen Books of the *Elements*.

Propositions I.5, 16, 26 are propositions by Thales, as mentioned by Proclus; the remainder of Book I, specifically the second part of Book I, I.28-47 is early Pythagorean, as it introduces the the Fifth Postulate and serves as an introduction to Book II.

Book VII, VII.1-19 is early Pythagorean, inspired by the Pythagorean arithmetization of music and based on the Pythagorean discovery of the Euclidean algorithm/anthyphairesis of natural numbers (Negrepontis, Farmaki, 2023). The remainder of Book VII is due to early Pythagoreans and most probably to Archytas, as it used exclusively for the Archytean Book VIII.

Book II has been reconstructed and restored to its original form aiming at the Pythagorean proof of incommensurability by means of the geometrical form in Negrepontis, Farmaki, 2025. Contributors to Books III, IV must include the Pythagoreans, Hippocrates, Theaetetus. Theaetetus theory of ratios of magnitudes (*) and Book X are due to Theaetetus, originally inspired from the Theodorus' incommensurabilities, as analyzed by Negrepontis and Protopapas, 2025 and present work. Book XIII, the natural continuation of Book X, and used by Plato in the *Timaeus*, is also the work of Theaetetus.

Archytas is unmistakably the creator of Book VIII with applications to music (as reported by Boethius, 1989) and arithmetical proofs of incommensurabilities, based on the second half of Book VII and continued for most part of Book IX.

The historically last part of the *Elements* is clearly the work by Eudoxus, Books V and VI on his theory of ratios of magnitudes, the precursor of the modern real numbers, and Books XI and XII, the precursor of the modern theory of integration.

> Book I (I.1-28 Thales, I.29-47 early Pythagoreans)
> Book VII (VII.1-19 early Pythagoreans)
> Book II (early Pythagoreans)
> Books III, IV (Pythagoreans, Hippocrates, Theaetetus)
> * Theaetetus' theory of proportion for magnitudes (Theaetetus)
> Book X (Theaetetus)
> Book XIII (Theaetetus)
> Books VII, 20-38, VIII, IX (Archytas)
> Books V, VI, XI, XII (Eudoxus)

| Table 2: Restoration of the historical order of the thirteen Books of the *Elements* | | | | | | | | | | | | |
|---|---|---|---|---|---|---|---|---|---|---|---|---|
| I | I | | | | | | | | | | | |
| II | | | II | | | | | | | | | |
| III | | | | III | | | | | | | | |
| IV | | | | | IV | | | | | | | |
| | | | | | | * | | | | | | |
| V | | | | | | | | | | | V | |
| VI | | | | | | | | | | | VI | |
| VII | | VII | | | | | | | | | | |
| VIII | | | | | | | | VIII | | | | |
| IX | | | | | | | | | IX | | | |
| X | | | | | | | X | | | | | |
| XI | | | | | | | | | | | | XI |
| XII | | | | | | | | | | | | XII |
| XIII | | | | | | | XIII | | | | | |

**Acknowledgment.** We express our warm thanks to Athanase Papadopoulos for his encouragement and extensive and essential suggestions that have led to a substantial improvement of our work.